\documentclass[11pt,reqno]{amsart}
\usepackage{euscript,amssymb,amsbsy}
\usepackage[arrow,matrix,curve]{xy}
\usepackage{a4wide}

\numberwithin{equation}{section}

\newenvironment{m-theorem}{%
\vskip4pt\refstepcounter{stff}\trivlist \itemindent 0pt
\item[\hskip\labelsep\bf Theorem \thestff]%
\it\ignorespaces}{\endtrivlist\vskip4pt}%
\newenvironment{m-proposition}{%
\vskip4pt\refstepcounter{stff}\trivlist \itemindent 0pt
\item[\hskip\labelsep\bf Proposition \thestff]%
\it\ignorespaces}{\endtrivlist\vskip4pt}%
\newenvironment{m-corollary}{%
\vskip4pt\refstepcounter{stff}\trivlist \itemindent 0pt
\item[\hskip\labelsep\bf Corollary \thestff]%
\it\ignorespaces}{\endtrivlist\vskip4pt}%
\newenvironment{m-lemma}{%
\vskip4pt\refstepcounter{stff}\trivlist \itemindent 0pt
\item[\hskip\labelsep\bf Lemma \thestff]%
\it\ignorespaces}{\endtrivlist\vskip4pt}%
\newenvironment{m-definition}{%
\vskip4pt\refstepcounter{stff}\trivlist \itemindent 0pt
\item[\hskip\labelsep\bf Definition \thestff]%
\ignorespaces}{\endtrivlist\vskip4pt}%
\newenvironment{m-notation}{%
\vskip4pt\refstepcounter{stff}\trivlist \itemindent 0pt
\item[\hskip\labelsep\bf Notation \thestff]%
\ignorespaces}{\endtrivlist\vskip4pt}%
\newenvironment{m-example}{%
\vskip4pt\refstepcounter{stff}\trivlist \itemindent 0pt
\item[\hskip\labelsep\bf Example \thestff]%
\ignorespaces}{\endtrivlist\vskip4pt}
\newenvironment{m-remark}{%
\vskip4pt\refstepcounter{stff}\trivlist \itemindent 0pt
\item[\hskip\labelsep\bf Remark \thestff]%
\ignorespaces}{\endtrivlist\vskip4pt}
\newenvironment{m-question}{%
\vskip4pt\refstepcounter{stff}\trivlist \itemindent 0pt
\item[\hskip\labelsep\bf Question.]%
\ignorespaces}{\endtrivlist\vskip4pt}%
\newenvironment{thm-nono}{
\vskip4pt\trivlist \itemindent 0pt
\item[\hskip\labelsep\bf Theorem.]%
\it\ignorespaces}{\endtrivlist\vskip4pt}%
\newenvironment{m-thank}{%
\vskip2pt\trivlist \itemindent 0pt
\item[\hskip\labelsep\it Acknowledgments]%
\ignorespaces}{\endtrivlist\vskip0pt}%
\newenvironment{m-proof}{%
\vskip4pt\trivlist \itemindent 0pt
\item[\hskip\labelsep\it Proof.]%
\ignorespaces}{\hfill$\Box$\endtrivlist\vskip4pt}%

\let\rar\rightarrow
\let\lar\longrightarrow

\let\hra\hookrightarrow
\let\mt\mapsto
\let\lmt\longmapsto

\font\tenmsa=msam10 %
\newcommand\hdashpiece{%
{\vrule height2.75pt depth-2.35pt width2.3pt \kern1.7pt}}%
\newcommand\hdashpieces{%
{\hdashpiece\hdashpiece\hdashpiece\hdashpiece}}%
\newcommand\dashto{\mathrel{%
\hdashpiece\hdashpiece\kern-0.4pt\hbox{\tenmsa K}}}%
\newcommand\dashar{\mathrel{%
\hdashpieces\kern-0.4pt\hbox{\tenmsa K}}}%

\let\euf\EuScript 
\let\cal\mathcal
\let\mbb\mathbb

\let\mfrak\mathfrak
 
\DeclareFontFamily{OT1}{rsfs}{}
\DeclareFontShape{OT1}{rsfs}{n}{it}{<->rsfs10}{}
\DeclareMathAlphabet{\crl}{OT1}{rsfs}{n}{it}

\let\ovl\overline

\let\tld\tilde

\let\nit\noindent
\let\disp\displaystyle
\let\srel\stackrel

\let\veps\varepsilon

\newcommand\lran[1]{{\langle #1\rangle}}

\newcommand\ort{\mathrel{{\vrule width4.0pt height0.4pt depth0pt
                \vrule width0.4pt height6.0pt depth0pt\,}}}
 
\newcommand\ev{{\rm ev}}
\newcommand\Aut{\mathop{\textrm{Aut}\kern1pt}}
\newcommand\cAut{\mathop{\mathcal{A}\kern-1pt\textit{ut}\kern1pt}}
\newcommand\End{\mathop{\rm{End}\kern1pt}}
\newcommand\cEnd{\mathop{\mathcal{E}\kern-1pt\textit{nd}\kern1pt}}
\newcommand\Hom{\mathop{\rm{Hom}\kern1pt}}
\newcommand\cHom{\mathop{\mathcal{H}\kern-1pt\textit{om}\kern1pt}}
\newcommand\Rad{\mathop{\rm{Rad}\kern1pt}}

\newcommand\Pic{\mathop{\rm Pic}\nolimits}
\newcommand\Spec{\mathop{\rm Spec}\nolimits}

\newcommand\Proj{\mathop{\rm Proj}\nolimits}

\newcommand\Ker{{\rm Ker}}

\newcommand\Img{{\rm Im}}
\newcommand\rk{{\rm rk}}

\let\si\sigma

\let\sm\setminus

\newcommand\bbY{{\mbb Y}}
\newcommand\bbk{\mbox{\rm I\kern-1.5pt k}}
\newcommand\sbbk{\hbox{\scriptsize I{\kern-.8pt}k}}
\newcommand\bone{{1\kern-0.57ex\rm l}}

\newcommand\bh{{\bf{h}}}
\newcommand\bk{{\bf{k}}}

\newcommand\eA{{\euf A}}
\newcommand\eB{{\euf B}}
\newcommand\eC{{\euf C}}

\newcommand\eF{{\euf F}}
\newcommand\eG{{\euf G}}
\newcommand\eI{{\euf I}}

\newcommand\cK{{\cal K}}

\newcommand\eL{{\euf L}}

\newcommand\eN{{\euf N}}

\newcommand\eO{{\euf O}}

\newcommand\cU{{\cal U}}
\newcommand\eU{{\euf U}}
\newcommand\cY{{\cal Y}}

\newcommand\eT{{\euf T}}
\newcommand\eW{{\euf W}}

\newcommand\pr{\mathop{\rm pr}\nolimits}
\newcommand\Flag{{\rm Fl}}
\newcommand\Grs{{\rm Gr}}
\newcommand\Sym{\mathop{\rm Sym}\nolimits}
\newcommand{\Ext}{\mathop{\rm Ext}\nolimits}

\newcommand\chr{{\rm char}}
\newcommand\codim{{\rm codim}}

\newcommand\Bl{{\rm Bl}}
\newcommand\Gl{{\rm Gl}}

\renewcommand\det{{\rm det}}

\let\ges\geqslant
\let\les\leqslant

\newcommand\ouset[3]{{\overset{#2}{\underset{#1}#3}}}

\newcommand\res{{\rm res}}
\newcommand{\cd}{\mathop{\rm cd}\nolimits}

\author{Mihai Halic}
\email{mihai.halic@gmail.com}

\begin{document}

\title[Split vector bundles along ample subvarieties]%
{Vector bundles on projective varieties whose restrictions to ample subvarieties split}

\begin{abstract}
We systematically study the splitting of vector bundles on a smooth, projective variety, whose restriction to the zero locus of a regular section of an ample vector bundle splits. 

First, we find ampleness and genericity conditions which ensure that the splitting of the vector bundle along the subvariety implies its global splitting. Second, we obtain a simple splitting criterion for vector bundles on the Grassmannian and on partial flag varieties.
\end{abstract}

\maketitle

\section*{Introduction}

We say that a vector bundle splits if it is isomorphic to a direct sum of line bundles.  
Horrocks proved in \cite{ho} his celebrated criterion for vector bundles on 
projective spaces, and his ideas gave rise to two main methods for proving the splitting 
of a vector bundle: either by imposing cohomological conditions, or by restricting them to 
hypersurfaces in the ambient space.

In this article we will follow the latter path. Although there are several splitting criteria 
obtained by restricting to divisors, it seems there are \emph{no similar results} for restrictions to 
higher co-dimensional subvarieties. Horrocks' result implies that a vector bundle on the 
projective space $\mbb P^d_\bk$---where $\bk$ is an algebraically closed field and 
$d\ges 3$---splits if an only if its restriction to a plane 
$\mbb P^2_\bk\subset\mbb P^d_\bk$ does so. Clearly, any plane is an ample subvariety 
of $\mbb P^d_\bk$. Our goal is to generalize this observation. Given a vector 
bundle $\crl V$ on a smooth, projective variety $X$, we ask under which assumptions 
the splitting of $\crl V$ along the zero locus of a regular section of an ample vector 
bundle $\eN$ on $X$ implies its global splitting. We investigate this issue from two points 
of view, each being interesting in its own right. 

First we prove that, if $\eN$ is \emph{sufficiently ample} compared to $\cEnd(\crl V)$, 
the splitting of $\crl V$ along the zero locus $Y_s$ of an \emph{arbitrary} regular section 
$s\in\Gamma(X,\eN)$ implies its global splitting. The proof requires the vanishing of various 
cohomology groups, and we carefully \emph{control} the amount of ampleness of $\eN$ 
necessary to achieve it. When the rank of $\eN$ is low compared to the dimension of $X$, 
our criteria take simple form, making them suited for concrete applications. 

Second, we avoid imposing ampleness on $\eN$; rather we focus 
on the \emph{genericity} of the section $s\in\Gamma(X,\eN)$. 
The splitting of $\crl V$ along the zero locus 
$Y_s\subset X$ implies its global splitting under the following hypotheses:
\begin{itemize}
\item[(i)] 
the ample vector bundle $\eN$ is globally generated, its rank is less than $\frac{1}{3}\dim X$ 
(cf. \eqref{eq:3nu}), and $s$ is \emph{very general} (in a precise sense);
\item[(ii)]
either the first cohomology group of any line bundle on $X$ vanishes, or the 
(finite dimensional) algebra of endomorphisms $\crl V\otimes\eO_{Y_s}$ is semi-simple.
\end{itemize}
If only (i) is fulfilled, then $\crl V$ is a successive extension of line bundles on $X$. 

In the last section we obtain a simple splitting criterion for vector bundles on the Grassmannian 
$\Grs(e;\bk^d)$ and on the partial flag variety $\Flag(d_1,\ldots,d_t;\bk^{d_{t+1}})$: 
a vector bundle splits if and only if, respectively, its restriction to an embedded $\Grs(2,\bk^4)$ 
and $\Flag(2,\ldots,2t;\bk^{2(t+1)})$ splits. 
We point out that, to our knowledge, currently there are \emph{no splitting criteria} 
(cohomological, uniformity, \textit{etc}.) for partial flag varieties. 

Let us elaborate on the results. Sufficient conditions which allow to extend the splitting of a vector bundle 
from a subvariety to the ambient space are given in 
the Proposition \ref{prop:formal}. It is the common root of the results obtained in this article and 
states \emph{without restrictions} that $\crl V$ splits on $X$ if and only if it does on 
the $m^{\rm th}$-order thickening $Y_{s,m}$ of the subvariety $Y_s$, of dimension at least one, 
for $m\gg0$. Care is taken to make this statement {\em effective}: 
we determine the order of the thickening of $Y_s$ in $X$ for which the splitting of $\crl V$ 
along $Y_{s,m}$ implies its global splitting. The proofs rely on the Buchsbaum-Eisenbud 
generic free resolutions \cite{be}, combined with  the vanishing theorems of Laytimi \cite{la} 
and Manivel \cite{ma}. For $Y$ smooth, we significantly improve the Proposition in two directions (cf. Theorem \ref{thm:formal+}): first, we state it \emph{intrinsically} for subvarieties $Y\subset X$ with ample normal bundle (thus justifying the term `ample subvariety' used in the title); second, we improve the previous bound. 

The Sections \ref{sct:2thm}, \ref{sct:gener} investigate the possibility of restricting $\crl V$ 
to the zero locus $Y_s$ itself, rather than to its thickening. The conclusions are contained in the 
Theorems \ref{thm:main1} and \ref{thm:generic-s2}: the former imposes \emph{ampleness conditions} for $\eN$ which are sufficient to restrict to zero loci of \emph{arbitrary} sections; the latter imposes \emph{genericity conditions} for $s$, and holds for $\eN$ ample, globally generated. 

Despite sharing a common root, the proofs of these results are \emph{very different} in nature. 
The cohomological criterion is based on \emph{effective cohomology vanishing} theorems; 
the genericity criterion is essentially a \emph{gluing} argument. 
We illustrate our results with concrete examples.

The last section discusses the necessity of the overall assumption that the vector bundle $\eN$, 
where we take sections, is ample. By analysing the case of the Grassmannian and of the partial 
flag varieties,  we conclude that the hypothesis can be weakened if $\eN$ is globally generated 
and one has enough control on the cohomology of the zero locus $Y_s$. The main results are the 
Theorems \ref{thm:grass} and \ref{thm:flag}, respectively, which massively 
simplify the complexity of the problem concerning the splitting of the vector bundles. 
The existing cohomological criteria \cite{ott,mal} for the Grasmannian involve a \emph{large} 
number of tests and, as far as we know, similar results are missing for flag varieties.


\section{The framework and the approach to the problem}
\label{sct:setup}

\begin{m-notation}\label{hyp}
Throughout the article $X$ stands for an irreducible, smooth, projective variety, defined over 
an algebraically closed field $\bk$ of characteristic zero. For a closed subscheme $S\subset X$, 
we denote by $\eI_S\subset\eO_X$ its sheaf of ideals. A {\it vector (resp. line) bundle} 
stands for a {\it locally free (resp. invertible) sheaf}. 

We consider two vector bundles $\crl V$ and $\eN$ on $X$, of rank $r$ and $\nu$ 
respectively, and assume that $\eN$ is ample. Let $\crl E:=\cEnd(\crl V)$ be 
the bundle of endomorphisms of $\crl V$, $\crl V_S:=\crl V\otimes\eO_S$, 
and similarly $\crl E_S:=\crl E\otimes\eO_S$. 

The zero locus of a section $s\in\Gamma(X,\eN)$ is the subscheme $Y_s$ of $X$ 
defined by the ideal sheaf $\eI_{Y_s}:={\rm Image}(\eN^\vee\srel{s}{\rar}\eO_X)$. 
We say that $s$ is a {\it regular section} if its zero locus $Y_s$ is a locally 
complete intersection in $X$. 
\end{m-notation}

Let $Y\subset X$ be the zero locus of a regular section of $\eN$, whose ideal sheaf is 
$\eI_Y$. For $m\ges 0$, the $m$-th order thickening $Y_m$ of $Y$ is the closed 
subscheme defined by $\eI_Y^{m+1}$; note that $Y_0=Y$ with 
this convention. The structure sheaves of two consecutive thickenings of $Y$ fit into 
the exact sequence 
\begin{equation}\label{eq:Ym}
0\to\Sym^m(\eN_Y^\vee)\cong\eI_Y^m/\eI_Y^{m+1}
\to\eO_{Y_{m}}\to\eO_{Y_{m-1}}\to 0.
\end{equation}

\begin{m-definition}\label{def:hat-xy}
The \emph{formal completion} of $X$ along $Y$ is defined as the direct limit 
$\varinjlim Y_m$, and it is denoted $\hat X_Y$. If no confusion is possible, 
we write $\hat X$. 
\end{m-definition}
When the ground field $\bk=\mbb C$, there are two kinds of thickenings and 
completions: one using germs of regular functions and another using germs of 
analytic functions. However, the generators of $\eI_Y^m$ are the same for all $m$, 
by Chow's theorem, so the exact sequence \eqref{eq:Ym} is valid in both cases. 

\begin{m-definition}\label{def:splitV}
We say that the vector bundle $\crl V$ \emph{splits} (or \emph{is split}) if there are $r$ 
line sub-bundles $\eL_1,\ldots,\eL_r\in\Pic(X)$ of $\crl V$ such that 
$\crl V =\mbox{$\ouset{i=1}{r}{\bigoplus}$}\,\eL_i$. Thus $\crl V$ splits if and 
only if there are pairwise non-isomorphic line sub-bundles $\eL_j$, $j\in J$, of $\crl V$ 
such that 
\begin{equation}
\label{eq:split}
\crl V =\mbox{$\ouset{j\in J}{}{\bigoplus}$}\,\eL_j\otimes\bk^{m_j}
\quad\text{with}\quad
\mbox{$\ouset{j\in J}{}{\sum}$}m_j=r.
\end{equation}
We call the vector sub-bundles $\crl V_j:=\eL_j\otimes\bk^{m_j}$, ${j\in J}$, 
the {\em isotypical components} of $\crl V$ corresponding to the splitting \eqref{eq:split}. 
\end{m-definition}

If $\mbox{$\ouset{j\in J}{}{\bigoplus}$}\,\eL_j\otimes\bk^{m_j}$ and 
$\mbox{$\ouset{j\,'\in J'}{}{\bigoplus}$}\,\eL'_{j\,'}\otimes\bk^{m\,'_{j\,'}}$
are two splittings of $\crl V$, then there is a bijective function 
$\si:J\to J'$ such that $\eL'_{\si(j)}\cong\eL_j$ and $m'_{\si(j)}=m_j$ 
for all $j\in J$. (See \cite[Theorem 1 and 2]{at}.) 

Unfortunately, the isotypical components are {\em not uniquely defined}, they depend 
on the choice of the splitting. Indeed, the global automorphisms of $\crl V$ send a 
splitting into another one. We define the following relation on the index set $J$: 
\begin{equation}\label{eq:order}
i\prec j\quad\Leftrightarrow\quad 
i\neq j\text{ and }\Gamma(X,\eL_i^{-1}\eL_j)\neq 0.
\end{equation}
It is straightforward to check that `$\prec$' is a partial order. The maximal elements with 
respect to $\prec$ have the property that the corresponding isotypical components are 
uniquely defined. 

\begin{m-lemma}\label{lm:max}
Let $M\subset J$ be the subset of maximal elements with respect to $\prec$. 
Then there is a natural, injective homomorphism of vector bundles 
\begin{equation}\label{eq:evm}
\ev_M:\underset{j\in M}{\mbox{$\bigoplus$}}
\eL_j\otimes\Gamma(X,\eL_j^{-1}\crl V)
\to \crl V.
\end{equation}
\end{m-lemma}

\begin{m-proof}
Clear, by the very definition.
\end{m-proof}

In this article we attempt to address the following: 
\begin{m-question}
Let $Y$ be the zero set of a regular section $s$ of $\eN$ and assume that $\crl V_{Y}$ splits. 
When does $\crl V$ split too?
\end{m-question}

The key to test the splitting of a vector bundle is the following elementary lemma, which allows to lift the splitting along a subvariety to the ambient space. 

\begin{m-lemma}\label{lm:eigv}
Let $S\subset X$  be a closed subscheme such that $\crl V_S$ splits and 
$$\res_S:\Gamma\bigl(X,\crl E\bigr)\to\Gamma\bigl(S,\crl E_S\bigr)$$ 
is surjective---in particular, $H^1(\eI_S\otimes\crl E)=0$. Then $\crl V$ splits. 
\end{m-lemma}

\begin{m-proof}
The hypothesis says that $\crl V_S\cong \ell_1\oplus\ldots\oplus\ell_r$, where 
$r:=\rk(\crl V)$ and $\ell_1,...,\ell_r\in\Pic(S)$. Take $\veps_1,...,\veps_r\in\bk$ 
pairwise distinct, and consider $\phi\in\End(V_S)$ given by multiplication by 
$\veps_\rho$ on $\ell_\rho$. Since $\res_S$ is surjective, $\phi$ extends to  
$\Phi\in\Gamma(X,\crl E)$. But the eigenvalues of $\Phi$ are independent of $x\in X$, 
so they are precisely $\veps_1,...,\veps_r$. Overall, we obtain $\Phi\in\End(\crl V)$ with 
$\rk(\crl V)$ distinct eigenvalues. Hence $\Phi_x$ is diagonalizable for all $x\in X$, 
and $\eL_\rho:=\Ker(\veps_\rho\bone-\Phi)$ are line bundles on $X$ such that 
$\crl V=\eL_1\oplus...\oplus\eL_r$. 
\end{m-proof}

\nit{\bf Strategy} \label{strategy}\quad 
To answer the question above we should show the surjectivity of the homomor\-phism $\res_Y:\Gamma(X,\crl E)\,{\to}\,\Gamma(Y,\crl E_Y)$. This is done in two stages:
\begin{enumerate}
\item Prove that $\Gamma(X,\crl E)\,{\to}\,\Gamma(Y_m,\crl E_{Y_m})$ is surjective, 
for $m\gg0$. Actually, we will determine \emph{effective lower bounds} for $m$ such that 
$H^1(\crl E\otimes\eI_Y^m)=0$. 

\item Prove that  $\Gamma(Y_m,\crl E_{Y_m})\,{\to}\,\Gamma(Y,\crl E_Y)$ is surjective, for all $m\ges 1$. 
\end{enumerate}

In this article, we will repeatedly use base change arguments: the property of a vector bundle to be split is unaffected by changing the (algebraically closed) ground field. The proposition below can be interpreted as the invariance of the Krull-Schmidt decomposition (see \cite{at}) under the change of the ground field. 

\begin{m-proposition}\label{prop:go-down}
Let $\bh,\bk$ be two algebraically closed fields of characteristic zero, such that 
$\bh\subset\bk$. Consider an irreducible projective scheme $X_\bh$ over $\bh$, 
and $\crl V_\bh$ a vector bundle on it. We define 
$X_\bk:=X_\bh\times_\bh\bk$ and $\crl V_\bk:=\crl V_\bh\times_{X_\bh}X_\bk$. 
Then $\crl V_\bh$ splits if and only if $\crl V_\bk$ splits.
\end{m-proposition}

\begin{m-proof}
We denote by $\Aut(\crl V_\bk)\subset\End(\crl V_\bk)$ the automorphism group of 
$\crl V_\bk$, and similarly for $\crl V_\bh$. 
Then $\Aut(\crl V_\bk)$ and $\Aut(\crl V_\bh)$ are linear algebraic groups 
and $\Aut(\crl V_\bk)=\Aut(\crl V_\bh)\otimes\bk$, by base change. 
Obviously, $\crl V_\bk$ splits if and only if the dimension of the maximal torus 
of $\Aut(\crl V_\bk)$ equals $\rk(\crl V_\bk)$, and the same holds for $\crl V_\bh$. 
But the dimension of the maximal torus is preserved under base change, and the 
conclusion follows.
\end{m-proof}


\section{The surjectivity of $\Gamma(X,\crl E)\to\Gamma(Y_m,\crl E_{Y_m})$ and effective lower bounds for $m$} \label{sct:x-xy}

Let $Y:=Y_s\subset X$ be the zero locus of a regular section $s\in\Gamma(X,\eN)$. 
Then its sheaf of ideals $\eI_Y\subset\eO_X$ admits the following well-known 
Koszul resolution:
\begin{equation}\label{eq:koszul}
0\lar 
\overset{\nu}{\hbox{$\bigwedge$}}\,\eN^\vee
\srel{s\ort}{\lar}
\overset{\nu-1}{\hbox{$\bigwedge$}}\,\eN^\vee
\lar\ldots\lar
\eN^\vee\srel{s\ort}{\lar}\eI_Y\lar 0.
\end{equation}
(Here $\ort$ stands for the contraction operation.) 
More generally,  locally free resolutions of the powers of $\eI_Y$ are constructed 
in \cite[Theorem 3.1]{be}. For any $m\ges 1$, we have the resolution 
\begin{equation}\label{eq:koszul-m}
0{\to} L^\nu_m(\eN^\vee){\to} L^{\nu-1}_m(\eN^\vee)
{\to}...{\to} L^j_m(\eN^\vee){\to}...{\to} 
\Sym^m(\eN^\vee)\srel{s^m\ort}{-\kern-1ex\lar}\eI_Y^m{\to} 0,
\end{equation}
where the vector bundles $L^j_m(\eN^\vee)$, $1\les j\les\nu$, are defined as follows:
\begin{equation}\label{eq:L}\kern-2ex
\begin{array}{rl}
L^j_m(\eN^\vee):=&\Ker\biggl(
\Sym^m(\eN^\vee)\otimes\overset{j-1}{\hbox{$\bigwedge$}}\,\eN^\vee 
\lar 
\Sym^{m+1}(\eN^\vee)\otimes\overset{j-2}{\hbox{$\bigwedge$}}\,\eN^\vee 
\biggr)
\\ 
=&
\Img\biggl(
\Sym^{m-1}(\eN^\vee)\otimes\overset{j}{\hbox{$\bigwedge$}}\,\eN^\vee 
\lar 
\Sym^{m}(\eN^\vee)\otimes\overset{j-1}{\hbox{$\bigwedge$}}\,\eN^\vee 
\biggr).
\end{array}
\end{equation} 
Actually $L^j_m(\eN^\vee)$ is a direct summand in both 
$\Sym^m(\eN^\vee)\otimes\overset{j-1}{\hbox{$\bigwedge$}}\,\eN^\vee$ and 
$\Sym^{m-1}(\eN^\vee)\otimes\overset{j}{\hbox{$\bigwedge$}}\,\eN^\vee$ 
because the homomorphisms which define $L^j_m(\eN^\vee)$ are 
$\cAut(\eN^\vee)$-invariant and the general linear group is linearly reductive. 
The long exact sequence \eqref{eq:koszul-m} breaks up into $\nu-1$ short exact 
sequences of the form 
\begin{equation}\label{eq:s}
0\to\euf S_{j+1}^{(m)}\to L^j_m(\eN^\vee)
\to\euf S_{j}^{(m)}
\to 0,\quad j=1,\ldots,\nu-1,
\end{equation}
with $\euf S_1^{(m)}=\eI_Y^m$ and 
$\euf S_{\nu}^{(m)}=\Sym^{m-1}(\euf N^\vee)\otimes\det(\eN^\vee)$. 

\begin{m-lemma}\label{lm:0}
Let $\euf F$ be a vector bundle of rank $f$ on $X$.\\[1ex]  
{\rm (i) (arbitrary $\nu$, lot of positivity for $\eN$)}\\ 
Let $m\ges 0$ be such that 
$
\Sym^{1+f}(\euf F^\vee)\otimes\det(\euf F)\otimes
\Sym^{m+\nu}(\eN)\otimes\det(\eN)^{-1}
$
is ample. Then holds 
$
H^t\bigl(X,\euf F\otimes\Sym^{m}(\eN^\vee)
\otimes\hbox{$\overset{j}\bigwedge$}\,\eN^\vee\bigr)=0,
\;\forall\,t<\dim X-\nu+j.
$\\ 
In particular, if 
$\Sym^{1+f}(\euf F^\vee)\otimes\det(\euf F)\otimes
\Sym^{1+\nu}(\eN)\otimes\det(\eN)^{-1}$ 
is ample and $\nu\les\dim X-2$, then  
$
H^j\bigl(X,\euf F\otimes\hbox{$\overset{j}\bigwedge$}\eN^\vee\bigr)=0, 
\forall\,j=1,\ldots,\nu.
$\\[1ex]
{\rm(ii) (low $\nu$, little positivity for $\eN$)}\\ 
If $\euf F^\vee\otimes\eN$ is ample and $\frac{(\nu+1)^2}{4}\les\dim X{-}f$, then 
$H^j\bigl(
\euf F\otimes\hbox{$\overset{j}\bigwedge$}\eN^\vee
\bigr)=0$, for all $j=1,\ldots,\nu$. 
\end{m-lemma}

\begin{m-proof}
(i) We consider the diagram 
$$
\xymatrix@C=2.5em@R=.5em{
&Z\ar[dd]^-{p}\ar[dr]\ar[dl]&
\\ 
\mbb P(\euf F)\ar[dr]^-{p_F}&&\mbb P(\eN^\vee)\ar[dl]_-{p_N}
\\ 
&X&
}
$$
where $Z:=\mbb P(\euf F)\times_X\mbb P(\eN^\vee)$ and 
$\eO_{p_F}(1)$, $\eO_{p_N}(1)$ stand for the relatively ample line bundles 
on $\mbb P(\eF)$, $\mbb P(\eN^\vee)$, respectively. 
Then the relative canonical bundle of $p_N$ satisfies 
$$\begin{array}{l}
\kappa_{p_N}
\!=
\det(\eN)\otimes\eO_{p_N}(-\nu), 
\text{ so}
\\[1ex] 
\Sym^m(\eN)
\!=
{(p_N)}_*\bigl(
\kappa_{p_N}\otimes\eO_{p_N}(m+\nu)\otimes\det(\eN)^{-1}
\bigr),\;\forall\,m\ges0.
\end{array}$$
Similar conclusion holds for the relative canonical bundle $\kappa_{p_F}$ of $p_F$. 

By hypothesis,  
$\eL:=\bigl(\eO_{p_F}(1+f)\otimes\det(\euf F)\bigr)
\boxtimes
\bigl(\eO_{p_N}(m+\nu)\otimes\det(\eN)^{-1}\bigr)$ 
on $Z$ is ample, and 
$\kappa_X\otimes\euf F^\vee\otimes\Sym^m(\eN)
=p_*(\kappa_Z\otimes\eL)$; the projection formula implies 
$$
H^i\bigl(X,
\kappa_X\otimes\euf F^\vee\otimes\Sym^m(\eN)\otimes
\hbox{$\overset{j}\bigwedge$}\eN\bigr)
=
H^i\bigl(
Z,\kappa_Z\otimes\eL\otimes\hbox{$\overset{j}\bigwedge$}\,p^*\eN
\bigr).
$$
On the right-hand-side, $p^*\eN$ is nef and $\eL$ is ample, so \cite{ma} 
implies that the cohomology group above vanishes for $i>\nu-j$. 
We conclude by applying the Serre duality on $X$. 

For the second claim: 
$H^j\big(\euf F^\vee\otimes\eN^\vee\otimes\hbox{$\overset{j-1}\bigwedge$}\eN^\vee\big)=0$, 
since $j<\dim X-\nu+j-1$ and $\hbox{$\overset{j}\bigwedge$}\eN$ 
is a direct summand in $\eN\otimes\hbox{$\overset{j-1}\bigwedge$}\eN$. 

\nit(ii) If $\euf F^\vee\otimes\eN$ is ample, then 
$\euf F^\vee\otimes\overset{j}{\bigwedge}\eN$ is ample too, being a direct summand 
of  $\euf F^\vee\otimes\eN^{\otimes j}$. 
Hence \cite[Theorem 2.1]{la} yields 
$H^{\dim X-j}(X,\kappa_X\otimes\euf F^\vee\otimes\overset{j}{\bigwedge}\eN)=0$, 
for all $j=1,\ldots,\nu$.
\end{m-proof}

\begin{m-proposition}\label{prop:cohom}
Let $\eN$ be an ample vector bundle of rank $\nu$ on $X$ and let 
$s\in\Gamma(X,\eN)$ be a regular section with zero locus $Y$. 
We consider an arbitrary, locally free sheaf $\euf F$ of rank $f$ on $X$. 
Then the following hold: 

\begin{itemize}
\item[(i)]  
Let $m_{\euf F}\ges 1$ be minimal such that 
$$
\Sym^{1+f}(\euf F^\vee)\otimes\det(\euf F)
\otimes\Sym^{m_{\eF}-1+\nu}(\eN)\otimes\det(\eN)^{-1}
$$
is ample. Then we have: 
$H^t(X,\euf F\otimes\euf I_Y^m)=0, \forall m\ges m_{\euf F},\forall t\les\dim Y{=}\dim X-\nu.$

In particular,  $H^j(X,\euf F)\to H^j(Y_m,\euf F)$ is an isomorphism, for $0\les j\les\dim Y-1$,  $m\gg0$. 

\item[(ii)]
Assume $f+\frac{(\nu+1)^2}{4}\les\dim X$, and $\euf F^\vee\otimes\eN$ is ample. 
Then $H^1(X,\euf F\otimes\eI_Y)=0$. 
\end{itemize}
\end{m-proposition}
 
\begin{m-proof}
(i) We tensor \eqref{eq:s} by $\euf F$. Since the middle term 
$\euf F\otimes L_m^j(\eN^\vee)$ is a direct summand in 
$\euf F\otimes\Sym^{m-1}(\eN^\vee)\otimes\overset{j}\bigwedge\eN^\vee\!$ 
for all $j$, the Lemma \ref{lm:0}(i) implies that for all $t\les\dim X-\nu$ holds 
$H^{t+j-1}\bigl(\euf F\otimes L^j_m(\eN^\vee)\bigr)=0$, hence 
$$
H^t(\euf F\otimes\eI_Y^m)\subset H^{t+1}(\euf F\otimes\euf S_{2}^{(m)})
\subset\ldots\subset
H^{t+\nu-2}(\euf F\otimes\euf S_{\nu-1}^{(m)})
\subset
H^{t+\nu-1}(\euf F\otimes L^{\nu}_{m}(\eN^\vee))=0.
$$
\nit(ii) The same argument as above, together with the Lemma \ref{lm:0}(ii) yields: \\ 
$H^1(X,\euf F\otimes\eI_Y)\subset H^2(X,\euf F\otimes\euf S_2^{(1)})
\subset\ldots\subset
H^\nu(X,\euf F\otimes\det\eN^\vee)=0.$
\end{m-proof}

\begin{m-corollary}\label{cor:ext}
Let $\bk=\mbb C$, and the situation be as in \ref{prop:cohom}, with $\dim Y\ges 2$. 
Consider an (analytic or Zariski) open neighbourhood $\cU$ of $Y$ in $X$, and let $\eA,\eC$ 
be two vector bundles on $X$. Then the following statements hold: 
\begin{enumerate}
\item[\rm(i)] Any extension of vector bundles on $\cU$, 
\\[1ex] 
$\null\hfill
0\to \eA\otimes\eO_\cU\to\eG_\cU\to\eC\otimes\eO_\cU\to 0,
\hfill\text{\rm(G)}$
\\[1ex] 
can be extended to an extension $0\to\eA\to\eB\to\eC\to 0$ of vector bundles 
on $X$, and $\eB$ is uniquely defined, up to isomorphism.

\item[\rm(ii)] Assume that the restriction to $\cU$ of three vector bundles $\eA,\eC,\eG$ 
on $X$ fit in the extension {\rm (G)}. Then $\eG$ is an extension of $\eC$ by $\eA$ on $X$. 
\end{enumerate}
\end{m-corollary}

\begin{m-proof}
(i)  Regardless whether the computations are done algebraically or analytically, the resolution \eqref{eq:koszul-m} and the Proposition \ref{prop:cohom} are both valid. 

The extension (G) corresponds to $\eta_\cU\in H^1(\cU,\eC^\vee\otimes\eA)$, 
and its restriction to $Y_m$ corresponds to the image 
$\eta_{m}\in H^1(Y_m,\eC^\vee\otimes\eA)$ of $\eta_\cU$. 
The Proposition \ref{prop:cohom}(i) implies that 
$H^1(X,\eC^\vee\otimes\eA)\to H^1(Y_m,\eC^\vee\otimes\eA)$, $m\gg0$, is an isomorphism, 
so $\eta_{m}$ uniquely lifts to $\eta\in H^1(X,\eC^\vee\otimes\eA)$; this 
defines the extension $0\to\eA\to\eB\to\eC\to 0$ on $X$. 
For the uniqueness part, notice that any two vector bundles $\eB$ and $\eB'$ are isomorphic 
along $Y_m$ and apply \ref{prop:cohom} again. 

\nit(ii) It is just a reformulation of the uniqueness statement above.
\end{m-proof}

\begin{m-lemma}\label{lm:split} 
Assume that $\crl V_{Y_{m_0}}$ splits, for some $m_0\ges 0$. Then  $\crl V$ splits as soon as either one of the following conditions is satisfied:

\nit{\phantom{ooo}\rm(i)} 
$\nu\les\dim X-1$, and 
$H^j\bigl(X,\crl E\otimes\Sym^{m_0}(\eN^\vee)\otimes \mbox{$\overset{j}\bigwedge$}\;\eN^\vee\bigr)=0$, for $j=1,\ldots,\nu$. 

\nit{or \rm(ii)} 
$\nu\les\dim X-2$, 
and $H^1\bigl(Y,\Sym^{m}(\eN_Y^\vee)\otimes\crl E_Y\bigr)=0$, for all $m\ges m_0+1$.
\end{m-lemma}

The twist of any vector bundle by a sufficiently ample line bundle satisfies the previous conditions. Horrocks' splitting criterion for $\mbb P^{\nu+2}$ is a particular case: just take $\eN:=\eO_{\mbb P^d}(1)^{\oplus\nu}$. The condition (ii) involves only $\crl E_Y$, which is a direct sum of line bundles; this simplifies the computations.\smallskip

\begin{m-proof}
(i) The hypothesis implies that $H^1\bigl(X,\crl E\otimes\eI_Y^{m_0+1}\bigr)=0$. 
(See the proof of \ref{prop:cohom}(i) above.) 

\nit(ii) The exact sequence \eqref{eq:Ym} implies that 
$\res^{Y_m}_{Y_{m-1}}{:}\,\Gamma(Y_m,\crl E){\to}\,\Gamma(Y_{m-1},\crl E)$ is surjective, for all $m\ges m_0+1$. Hence $\Gamma\bigl(X,\crl E\bigr)\to\Gamma\bigl(Y,\crl E_{Y_{m_0}}\bigr)$ is surjective too, by \ref{prop:cohom}(i). 
\end{m-proof}

The following general statement holds for arbitrary vector bundles on projective varieties. 
It is the common root of the subsequent results in this article. 

\begin{m-proposition}\label{prop:formal} 
Let $\eN$ be an ample vector bundle of rank $\nu\les\dim X-1$, and let $Y$ be the zero locus of a regular section of $\eN$. 

\nit{\rm (i)} Let $m_{\crl V}\ges0$ be minimal such that 
$\Sym^{1+r^2}(\crl E)\otimes\Sym^{m_{\crl V}+\nu}(\eN)\otimes\det(\eN)^{-1}$ is ample. 
Then $\crl V$ splits if and only if $\crl V_{Y_m}$ does, for some $m\ges m_{\crl V}$. 

In particular, $\crl V$ is split if and only if its restriction to $\hat X_Y$ splits. 

\nit{\rm (ii)} Assume that the ground field is $\bk=\mbb C$.  
Then $\crl V$ splits if and only if there is an open {\em analytic neighbourhood} $\cal U$ 
of $Y$ such that $\crl V\otimes\eO_{\cal U}$ splits. 
\end{m-proposition}

\begin{m-proof}
Apply the Proposition \ref{prop:cohom} and the Lemma \ref{lm:eigv}. 
\end{m-proof}

When the subvariety $Y\subset X$ is smooth, the bounds in previous statement can be improved, by using the bootstrapping argument \ref{lm:split}(ii) combined with the Kodaira vanishing theorem. Notably, the statement below is \emph{intrinsic} to $Y$, \emph{does not involve} the vector bundle $\eN$, thus it justifies the term of `ample subvariety' used in the title.

\begin{m-theorem}\label{thm:formal+}
Let $Y\subset X$ be a \emph{smooth}, irreducible subvariety and denote by $\eN_Y$ its normal bundle. We assume that the following hypotheses are satisfied:
\begin{enumerate}
\item[\rm(i)] 
\begin{enumerate}
\item[\rm(a)] 
$\nu:=\codim_X Y\les\dim X-2$; 
\item[\rm(b)] 
$\eN_Y$ and $\,\crl E_Y\otimes\Sym^{m_0+1+\nu}(\eN_Y) \otimes{(\det\,\eN_Y)}^{-1}$ are ample. 
\end{enumerate}
\item[\rm(ii)]
The cohomological dimension of the complement is $\cd(X\sm Y)\les\dim X-2$. \\ 
The inequality holds if $Y$ is the zero locus of a regular section in an ample vector bundle $\eN$ on $X$, of rank $\nu\les\dim X-1$. 
\end{enumerate}
Then $\crl V$ splits if and only if $\crl V_{Y_{m_0}}$ does.
\end{m-theorem}

\begin{m-proof}
We prove that (a) and (b) imply $H^1\bigl(Y,\Sym^{m}(\eN_Y^\vee)\otimes\crl E_Y\bigr)=0$, for all $m\ges m_0+1$, so $\disp\Gamma(\hat X_Y,\crl E_{\hat X_Y})=\varprojlim_m\Gamma(Y_m,\crl E)\to \Gamma(Y,\crl E_Y)$ is surjective (cf. \eqref{eq:Ym}). The projection formula yields: 
$$
H^1\bigl(Y,\Sym^m(\eN_Y^\vee)\otimes\crl E_Y\bigr)
\cong
H^\nu\bigl(
\mbb P(\eN_Y^\vee),
{\bigl(\eO_q(m+\nu)\otimes{(q^*\det\,\eN_Y)}^{-1}
\otimes q^*\crl E_Y\bigr)}^\vee
\bigr).
$$
For $\crl V_Y\cong\ouset{j=1}{r}{\oplus}\ell_j$, we have 
$\crl E_Y\cong\ouset{i,j}{}{\oplus}\,\ell_j\ell_i^{-1}$, hence 
$$
\begin{array}{r}
H^\nu\bigl(
\mbb P(\eN_Y^\vee),
{\bigl(\eO_q(m+\nu)\otimes{(q^*\det\,\eN_Y)}^{-1}
\otimes q^*\crl E_Y\bigr)}^\vee
\bigr)
\\[1ex]
=
\underset{\ell}{\bigoplus}\; 
H^\nu\bigl(
\mbb P(\eN_Y^\vee),
{\bigl(\eO_q(m+\nu)\otimes{(q^*\det\,\eN_Y)}^{-1}
\otimes q^*\ell\bigr)}^{-1}
\bigr).
\end{array}
$$
where $\ell$ runs over the direct summands of $\crl E_Y$. The terms on the right hand side vanish, because $\eO_q(m+\nu)\otimes{(q^*\det\,\eN_Y)}^{-1}\otimes q^*\ell$ is ample for $m\ges m_0+1$. 

The assumption (ii) implies that $\Gamma(X,\crl E)\to\Gamma(\hat X_Y,\crl E_{\hat X_Y})$ is surjective, by \cite[Theorem III.3.4(b)]{ha}, so overall $\Gamma(X,\crl E)\to\Gamma(Y_{m_0},\crl E)$ is surjective too. 

Finally, the claim about the condition (ii) of the theorem follows from \cite[Theorem III.3.4(b)]{ha} and the Proposition \ref{prop:cohom}(i).
\end{m-proof}

\nit In the subsequent sections we will obtain sufficient conditions for these criteria.


\section{The surjectivity of $\Gamma(Y_m,\crl E_{Y_m})\,{\to}\,\Gamma(Y,\crl E_{Y})$ and $1^{\text{st}}$ criterion:\break ampleness conditions for $\eN$} \label{sct:2thm} 

In this section we consider {\em arbitrary} regular sections of $\eN$ such that $\crl V$ splits along their zero locus, and we impose {\em sufficient ampleness} on $\eN$ in order to deduce the global splitting of $\crl V$. 

\begin{m-theorem}\label{thm:main1} 
The implication $\,[\,\crl V_Y\text{ splits }\Rightarrow\crl V\text{ splits}\,]\,$ holds in any of the situations described below. 
\begin{itemize}
\item[(a)] 
Assume that $\nu=\rk(\eN)\les\dim X-2$, and let $s\in\Gamma(X,\eN)$ be an \emph{arbitrary} regular section with zero locus $Y$. 

\begin{enumerate}
\item[(i)]
$\Sym^{1+r^2}(\crl E)\otimes\Sym^{1+\nu}(\eN)\otimes\det(\eN)^{-1}$ is ample; 

\item[(ii)]
$\frac{(\nu+1)^2}{4}\les\dim X-r^2$ and $\crl E\otimes\eN$ is ample;

\item[(iii)]
$Y$ is smooth, $\nu\les\frac{\dim X-1}{2}$, and $\eN=\eG\otimes\eA$ with $\eG$ a globally generated vector bundle of rank $\nu$, and $\eA$ an ample line bundle such that $\crl E_Y\otimes\eA_Y$ is ample.
\end{enumerate}

\item[(b)] 
$Y\subset X$ is a \emph{smooth} subvariety with the following properties: 
\begin{enumerate}
\item[\rm(i)] 
$\nu:=\codim_X Y\les\dim X-2$ and $\cd(X\sm Y)\les\dim X-2$; 
\item[\rm(ii)] 
$\eN_Y$ and $\,\crl E_Y\otimes\Sym^{1+\nu}(\eN_Y) \otimes{(\det\,\eN_Y)}^{-1}$ are ample.  
\end{enumerate}
\end{itemize}
\end{m-theorem}


\begin{m-proof} 
(a)(i) The claim follows from the Lemma \ref{lm:0}(i) and \ref{lm:split}(i).

\nit(ii) The Propositions \ref{prop:cohom}(ii) implies that $H^1(X,\crl E\otimes\eI_Y)=0$. 

\nit(iii) 
For all the direct summands $\ell$ of $\crl E_Y$, the line bundle $\ell\otimes\eA_Y$ is ample. Then 
$$
H^1(Y,\Sym^m(\eN_Y^\vee)\otimes\crl E_Y)^\vee
=
\underset{\ell}{\mbox{$\bigoplus$}}\;
H^{\dim Y-1}(\,Y,\kappa_Y\otimes\kern-1.25ex
\underbrace{\Sym^m(\eG_Y)}_{\text{\scriptsize globally generated}}\kern-1.25ex\otimes\,
\underbrace{(\ell\otimes\eA_Y^m)}_{\text{\scriptsize ample}}\,),
$$
vanishes, for all $m\ges1$, by \cite[Theorem 2.4]{la}. The Lemma \ref{lm:split}(ii) yields the conclusion.

\nit(b) This is the Theorem \ref{thm:formal+}, for $m_0=0$. 
\end{m-proof}

In some cases one wishes to prove the triviality of a vector bundle (cf. \cite{bs}). 

\begin{m-corollary}\label{cor:triv}
Assume that $\nu\les\frac{\dim X-1}{2}$ and $\eN=\eG\otimes\eA$, with $\eG$ 
globally generated and $\eA$ ample. If $\crl V$ is trivializable along the zero locus of 
a regular section in $\eN$, and this zero locus is smooth, then $\crl V$ is trivializable 
on $X$. 
\end{m-corollary}

\begin{m-proof}
Indeed, in this case $\crl E_Y\cong\eO_Y^{\oplus r^2}$.
\end{m-proof}

With the Theorem \ref{thm:main1}(iii) one can create a host of examples. To check the splitting of a vector bundle $\crl V$, one should proceed as follows: find (according to the case) a low rank, globally generated vector bundle $\eG$ on $X$, and an ample line bundle $\eA$ such that $\crl E\otimes\eA$ is ample; then restrict $\crl V$ to the zero locus of a section in $\eG\otimes\eA$. 

\begin{m-example}\label{expl:grass} 
(Compare with \ref{expl:2n}(i)) 
Let $X\srel{\iota}{\hra}\Grs(n;\mbb C^{n+\nu})$ be a $c$-codimen\-sional subvariety, 
with $n\ges3,\nu\ges2, c\les(n-2)\nu-1$; let $\eG$ be the universal quotient bundle 
and $\eO_X(1):=\det(\eG)_X$. The following statements hold:
\begin{enumerate}
\item[(i)] 
If $\crl E(a)$ is ample for some $a\ges 1$, and $\crl V$ splits along the (smooth) 
zero locus of an arbitrary regular section in $\eN:=\eG_X(a)$, then $\crl V$ splits on $X$. 
\item[(ii)]
If the restriction of $\crl V$ to the (smooth) zero locus of a regular section 
in $\eG_X(1)$ is trivializable, then $\crl V\cong\eO_X^{\oplus r}$. 
\end{enumerate}
Similar statements hold for the other (isotropic) Grassmannians, too.
\end{m-example}


\section{A gluing procedure and the $2^{\text{nd}}$ criterion:\break splitting along zero loci of generic sections of $\eN$} \label{sct:gener}

Now we change our viewpoint. Instead of imposing ampleness on $\eN$, we 
prove that the splitting of a vector bundle along the zero locus of a {\em very general} 
section of a globally generated ample vector bundle implies its global splitting. 
Throughout this section we assume that $\eN$ is {\em globally generated}, and furthermore: 
\begin{equation}\label{eq:3nu}
\null\kern-2.25ex
\begin{array}{cl}
& \kern-4ex
\nu\les\min\bigl\{
\frac{\dim X-3}{2},\frac{\dim X-1}{3}
\bigr\} 
\;\text{that is}\;
\left\{\begin{array}{ll}
\nu=1&\text{ for }\dim X=5,6,
\\
\nu\les\frac{\dim X-1}{3}&\text{ for }\dim X\ges 7,
\end{array}\right.
\\[3ex] 
\text{or}
&
\nu=1,\,\dim X=4,\text{ and }\kappa_X\otimes\eN^2\text{ is globally generated,} 
\\
&
\text{where $\kappa_X$ stands for the canonical bundle of $X$. }
\end{array}\kern-1.75ex
\end{equation}
Our goal is to prove that the splitting of $\crl V$ along the geometric generic section 
of $\eN$ implies its global splitting. 
The proof uses base change arguments, so we start with general considerations. 
The variety $X$ and the vector bundles $\eN,\crl V$ are defined by equations involving 
finitely many coefficients in $\bk$. After adjoining them to $\mbb Q$, we obtain a field 
extension of finite type $\mbb Q\hra\bk_0$. In particular, $\bk_0$ is countable, so we 
can realize it as a sub-field of $\mbb C$.
$$
\bk_0\hra\bk
\;\srel{\bk\text{ \scriptsize alg. closed}}{\Longrightarrow}\; 
\bar\bk_0\hra\bk
\qquad\text{and}\qquad 
\bk_0\hra\mbb C
\;\srel{\mbb C\text{ \scriptsize alg. closed}}{\Longrightarrow}\; 
\bar\bk_0\hra\mbb C.
$$
After replacing $\bk_0$ by $\bar\bk_0$, we find a countable, algebraically closed field 
$\bk_0$, which is simultaneously a sub-field of $\bk$ and of $\mbb C$, such that 
$X,\eN,\crl V$ are defined over $\bk_0$. In this situation we have the Cartesian, 
base change diagram 
$$
\xymatrix@R=1.5em{
X=X_\bk\ar[d]\ar[r]^-b
&X_0:=X_{\bk_0}\ar[d]
\\ 
\Spec(\bk)\ar[r]&\Spec(\bk_0)
}
$$
and there are vector bundles $\eN_{0}, \crl V_{0}$ on $X_{0}$ such that 
$\eN=\eN_{0}\times_{\bk_0}\bk$
and also 
$\crl V=\crl V_{0}\times_{\bk_0}\bk$; we denote $\crl E_0:=\cEnd(\crl V_0)$. 
By base change, $\eN_0$ on $X_0$ is globally generated too. 
Let $\mbb P^N_\bk:=\mbb P(\Gamma(X,\eN))=\Proj\bigl(
\,\Sym_\bk^\bullet(\Gamma(X,\eN)^\vee)
\bigr)$, and similarly for $\bk_0$, and we consider the trace morphism 
$$
\begin{array}{r}
\mbb P^N_\bk
\lar
\mbb P^N_{\bk_0},
\quad 
\mfrak p
\lmt
\mfrak p\cap \Sym_{\bk_0}^\bullet
\bigl(\Gamma(X_0,\eN_0)^\vee\bigr).
\end{array}
$$
The sheaf $\cK$ defined by 
\begin{equation}\label{eq:N}
0\to\cK:=\Ker(\eta)\to\Gamma(X,\eN)\otimes\eO_X\srel{\eta}{\to}\eN\rar0
\end{equation}
is locally free, and the incidence variety 
$\cY:=\{([s],x)\mid s(x)=0\}\subset\mbb P^N_\bk\times X$  is naturally isomorphic to 
the projective bundle $\mbb P(\cK)$ over $X$. We denote by $\pi$ and $q$ respectively 
the projections of $\cY$ onto $\mbb P^N_\bk$ and $X$. For any open subset $S$ of 
$\mbb P^N_\bk$, we let $\cY_S:=\pi^{-1}(S)$. If the ground field $\bk=\mbb C$, 
we will consider open subsets of $\mbb P^N_{\mbb C}$ in the {\em analytic} topology. 
Henceforth we use this notation. 

\begin{m-definition}\label{def:gg-section} 
Let $\bbk$ be the quotient field of $\mbb P^N_\bk$, and $\bar{\bbk}$ its algebraic closure. 
The {\em geometric generic section} $\mbb Y$ of $\eN$ is defined 
by the Cartesian diagram: 
$$
\xymatrix@R=1.5em{
\mbb Y:=\cY_{\bar{\sbbk}}\ar[d]\ar[r]^-\psi&\cY\ar[d]^-\pi
\\
\Spec(\bar{\bbk})\ar[r]&\mbb P^N_\bk.
}
$$
\end{m-definition}

The next lemma shows that the assumption that $\crl V$ splits on $\mbb Y$ is \emph{a priori weaker} than to say that $\crl V$ splits on $\cY\times_{\mbb P^N_{\bk}}\Spec(\bbk)$, that is on the generic section of $\eN$. The upshot of this section is to prove that the ampleness of $\eN$ actually forces $\crl V$ to split on $X$. We believe that this fact is indeed unexpected.

\begin{m-lemma}\label{lm:galois}
Assume that the restriction of $\crl V$ to the geometric generic section of $\eN$ splits. 
Then there is a non-empty Zariski open subset $S$ of $\mbb P^N_\bk$, and a finite cover 
$S'\to S$ such that $q^*\crl V\times_SS'$ splits on $\cY_{S'}$, and $Y_s$ is smooth 
for all $s\in S$. 
If $\bk=\mbb C$, there is an open {\em ball} $B\subset\mbb P^N_{\mbb C}$ 
with the previous two properties. 
\end{m-lemma}

\begin{m-proof}
Let $(q^*\crl V)_\bbY$ be the pull-back of $q^*\crl V$ to $\bbY$; 
there are $\ell'_1,\ldots,\ell'_r\in\Pic(\bbY)$ such that 
$(q^*\crl V)_\bbY=\ell'_1\oplus\ldots\oplus\ell'_r$. 
Since $\ell'_1,\ldots,\ell'_r$ are defined over an intermediate field 
$\bbk\hra\bbk'\hra\bar{\bbk}$ finitely generated and algebraic over $\bbk$, there is an open affine 
$S\subset\mbb P(\Gamma(X,\eN))$, an affine variety $S'$, and a finite morphism 
$S'\srel{\si}{\to} S$ such that $\ell'_1,\ldots,\ell'_r$ are defined over $\bk[S']$ and 
$(q^*\crl V)_{S'}$ splits on $\cY_{S'}$. 
After shrinking $S$ further, $Y_s$ is smooth for all $s\in S$, by Bertini's theorem.

If $\bk=\mbb C$, there are open balls $B'\subset S'$ and $B\subset S$ such that 
$\si:B'\to B$ is an isomorphism. Then the splitting of 
${(q^*\crl V)}_{B'}$ descends to ${(q^*\crl V)}_{B}$ on $\cY_{B}$.
\end{m-proof}

Henceforth we assume $\bk=\mbb C$ and consider an open ball $B\subset\mbb P\bigl(\Gamma(X,\eN)\bigr)$ 
as above. We choose an isotypical decomposition: 
\begin{equation}\label{eq:iso-qV}
{(q^*\crl V)}_{B}=\underset{j\in J}{\bigoplus}\,\ell_j\otimes\mbb C^{m_j}, 
\text{ with }\ell_j\in\Pic(\cY_{B}) \text{ pairwise non-isomorphic}. 
\end{equation}
For $(s,t)\in B\times B$, the intersection $Y_{st}:=Y_s\cap Y_t$ is the zero locus 
of $(s,t)\in\Gamma(X,\eN^{\oplus 2})$. Since $\eN$ is globally generated, 
$Y_s$ and $Y_t$ meet transversally for $(s,t)$ in an open, dense subset $(B\times B)^\circ$. 
We consider the diagram: 
\begin{equation}\label{eq:pic}
\xymatrix@R=.5em@C=4em{
&\Pic(Y_s)\ar[dr]^-{\res^{Y_s}_{Y_{st}}}&
\\ 
\Pic(X)\ar[ur]^-{\res^X_{Y_s}}
\ar[dr]_-{\res^X_{Y_t}}\ar[rr]|{\;\res^X_{Y_{st}}\;}
&&
\Pic(Y_{st}).
\\ 
&\Pic(Y_t)\ar[ur]_-{\res^{Y_t}_{Y_{st}}}&
}
\end{equation}
Assume that $\dim X\ges 5$. Then the Lefschetz-Sommese theorem \cite{so} implies that 
all the arrows are isomorphisms, for all $(s,t)\in (B\times B)^\circ$. 

Now assume that $\dim X=4$. Since 
$\kappa_{Y_s}\otimes\eN=(\kappa_X\otimes\eN^2)\otimes\eO_{Y_s}$ 
is globally generated,  the Noether-Lefschetz theorem \cite{rs} implies that the \eqref{eq:pic} 
consists of isomorphisms for a dense subset of $(B\times B)^\circ$.

\begin{m-lemma}\label{lm:pic}
The pull-back $\Pic(X)\srel{q^*}{\to}\Pic(\cY_{B})$ is an isomorphism, so 
$$
{(q^*\crl V)}_B\cong
q^*\Bigl(\underset{j\in J}{\mbox{$\bigoplus$}}\eL_j^{\oplus m_j}\Bigr)\otimes\eO_{\cY_{B}}, 
\text{ with }\eL_j\in\Pic(X).
$$
\end{m-lemma}

\begin{m-proof} 
Fix $o\in B$. 
The composition $\Pic(X)\srel{q^*}{\to}\Pic(\cY_B)\srel{\res_{Y_o}}{\to}\Pic(Y_o)$ 
is bijective, so $q^*$ is injective. 
For the surjectivity, take $\ell\in\Pic(\cY_B)$. If $\ell_{Y_o}\cong\eO_{Y_o}$, then 
$$
\{s\in B\mid \ell_{Y_s}\not\cong\eO_{Y_s}\}=\{s\in S\mid h^0(\ell_{Y_s})=0\}
$$
is open, by semi-continuity, so $\{s\in B\mid\ell_{Y_s}\cong\eO_{Y_s}\}$ is closed. 
On the other hand, by restricting to $Y_{os}$, the previous discussion implies that this 
set is dense; thus it is the whole $B$. 
It follows that $\ell\cong\pi^*\bar\ell$, with $\bar\ell\in\Pic(B)$, so $\ell\cong\eO$. 
If $\ell\in\Pic(\cY)$ is arbitrary, take $\eL\in\Pic(X)$ such that 
$\ell_{Y_o}\cong\eL_{Y_o}$, so ${(q^*\eL^{-1})\ell|}_{Y_o}$ is trivial. 
\end{m-proof}

For all $s\in B$, let $M_s\subset J$ 
be the subset of maximal elements with respect to \eqref{eq:order}, corresponding to the 
splitting of $\crl V\otimes \eO_{Y_s}$. By semi-continuity, for any $s\in B$, 
there is a neighbourhood $B_s\subset B$ of $s$ such that $M_s\subset M_{s'}$ for all 
$s'\in B_s$. Thus there is a largest subset $M\subset J$, and an open subset $B'\subset B$ 
such that $M=M_s$ for all $s\in B'$.  

\begin{m-lemma}\label{lm:ind-start}
Assume that $\bk=\mbb C$ and \eqref{eq:3nu} is satisfied. 
Furthermore, let $B\subset\mbb P^N_{\mbb C}$ be a ball such that $Y_s$ is 
smooth for all $s\in B$, ${(q^*\crl V)}_{B}$ splits over $\cY_B$, and the set of 
maximal elements $M\subset J$ with respect to $\prec$ is the same for all $s\in B$. 

We consider the (analytic) open subset $\cU:=q(\cY_B)\subset X$. 
Then there is an injective homomorphism of vector bundles  
$
{\bigl(
\underset{\mu\in M}{\mbox{$\bigoplus$}}\,
\eL_\mu^{\oplus m_\mu}
\bigr)}\otimes\eO_\cU
\to\crl V\otimes\eO_\cU
$ 
whose restriction to $Y_s$ is the natural evaluation \eqref{eq:evm}, for all $s\in B$. 
\end{m-lemma}

\begin{m-proof}
The restriction to $Y_s$ of  
$
\ev:\kern-.5ex
\underset{\mu\in M}{\mbox{$\bigoplus$}}\,
q^*\eL_\mu\otimes\pi^*\pi_*q^*(\eL_\mu^{-1}\otimes\crl V)_B
\to (q^*\crl V)_B
$
is the homomorphism \eqref{eq:evm}, for all $s\in B$. The maximality of $\mu\in M$ 
implies that $\pi_*q^*(\eL_\mu^{-1}\otimes\crl V)\cong\eO_B^{\oplus m_\mu}$ 
and $\ev$ is pointwise injective. We prove that, after suitable 
choices of bases in $\pi_*q^*(\eL_\mu^{-1}\otimes\crl V)$, $\mu\in M$, the 
homomorphism $\ev$ descends to $\cU$. We will deal with each index separately, 
the overall basis being the direct sum of the individual ones. 

Consider $\mu\in M$, and a base point $o\in B$. 
Then $\crl V':=\eL_\mu^{-1}\otimes\crl V$ has the following properties:
\begin{enumerate}
\item[--] 
$(q^*\crl V')_B\cong\eO_{\cY_{B}}^{\oplus m}\oplus
\Big(\,
\underset{j\in J\sm\{\mu\}}{\bigoplus}q^*{(\eL_\mu^{-1}\eL_j)}^{\oplus m_j}_B$
\Big). 
\item[--] 
$\pi_*(q^*\crl V')_B\cong\eO_{B}^{\oplus m}$; 
we choose an isomorphism $\alpha_B$ between them.
\item[--] 
$\pi^*\pi_*(q^*\crl V')_B\to (q^*\crl V')_B$ is pointwise injective; 
let $\eT\subset (q^*\crl V')_B$ be its image.
\end{enumerate}
\nit We choose a complement 
$\eW\cong\underset{j\in J\sm\{\mu\}}{\bigoplus}
q^*{(\eL_\mu^{-1}\eL_j)}^{\oplus m_j}$ of $\eT$ in $(q^*\crl V')_B$, that is 
\begin{equation}\label{eq:tw}
{(q^*\crl V')}_B=\eT\oplus\eW.
\end{equation}
The isomorphism $\alpha_B$ above determines the pointwise injective homomorphism 
$$
\alpha:\eO_{\cY_{B}}^{\oplus m}\to (q^*\crl V')_B=\eT\oplus\eW
$$
whose second component vanishes, as $\Gamma(\cY_{B},\eW)=0$. 
Let $\beta:(q^*\crl V')_B\to\eO_{\cY_{B}}^{\oplus m}$ 
be the left inverse of $\alpha$ with respect to the splitting \eqref{eq:tw}, 
and note $\alpha\circ\beta|_\eT=\bone_\eT$. 

\nit{\it Claim}\quad After a suitable change of coordinates 
in $\eO_{\cY_{B}}^{\oplus m}$, the homomorphisms $\alpha$ descends to 
$\cU=q(\cY_{B})\subset X$. Indeed, for any $s\in B$, we consider the diagram 
$$
\xymatrix@C=3em@R=1.5em{
\kern1ex\eO_{Y_{os}}^{\oplus m}\ar[r]^-{\alpha_o}\ar@{..>}[d]_-{a_s}
&
\crl V'_{Y_{os}}\ar@{=}[d]
&
\\ 
\kern1ex\eO_{Y_{os}}^{\oplus m}\ar[r]^-{\alpha_s}
&
\crl V'_{Y_{os}}\ar@<-2pt>`d[l] `[l]|-{\beta_s}
&
\text{with}\;
a_s:=\beta_s\circ\alpha_o\in\End(\mbb C^m).
}
$$
Similarly, we let $a'_s:=\beta_o\circ\alpha_s$. Then holds 
$a'_sa_s=\beta_o\alpha_s\beta_s\alpha_o=\beta_o\alpha_o=\bone$ 
(for the second equality notice 
$\Img\bigl({\alpha_o|}_{Y_{os}}\bigr)
=\eT_{Y_{os}}=\Img\bigl({\alpha_s|}_{Y_{os}}\bigr)$, 
and $\alpha_s\beta_s|_\eT=\bone$), and similarly $a_sa'_s=\bone$. 
Thus $a_s\in\Gl(m;\mbb C)$ for all $s\in B$, and the new trivialization 
$\tilde\alpha:=\alpha\circ a$ of $\eT$ satisfies  
\begin{equation}\label{eq:os}
\tilde\alpha_s=\tilde\alpha_o\text{ along }Y_{os},
\quad\forall\,s\in B, 
\end{equation}
because 
$\tilde\alpha_s|_{Y_{os}}{=}(\alpha_s\beta_s)\alpha_o|_{Y_{os}}
{=}\alpha_o|_{Y_{os}}{=}\tilde\alpha_o|_{Y_{os}}$. 
Now we observe that, for all $s,t\in B$, the trivializations of $\eT_{Y_{st}}$ induced by 
$\tilde{\alpha}$ from $Y_s$ and $Y_t$, coincide. Equivalently, the 
following diagram commutes: 
\begin{equation}\label{eq:st}
\xymatrix@R=1.25em@C=2.5em{
\eO_{Y_{st}}^{\oplus m}\ar[r]^-{\tilde\alpha_s}_-\cong\ar@{=}[d]
&
\eT_{Y_{st}}\ar@<-.5ex>@{=}[d]
\kern-1.5ex\ar@{}[r]|-{\mbox{$\subset$}}
&
\kern-1.5ex
\crl V'_{Y_{st}}\ar@<-2ex>@{=}[d]
&
\ar@{}[d]^-{\mbox{$
\Leftrightarrow\;
{\tilde\alpha_t^{-1}\circ\tilde\alpha_s|}_{Y_{st}}=\bone\in\Gl(r;\mbb C).
$}}
\\ 
\eO_{Y_{st}}^{\oplus m}\ar[r]^-{\tilde\alpha_t}_-\cong
&
\eT_{Y_{st}}
\kern-1.5ex\ar@{}[r]|-{\mbox{$\subset$}}
&
\kern-1.5ex
\crl V'_{Y_{st}}
&
}
\end{equation}
Indeed, the triple intersection $Y_{ost}$--- the zero locus of 
$(o,s,t)\in\Gamma(X,\eN^{\oplus 3})$---is a non-empty, connected subscheme of $X$, 
as $\dim X-3\nu\ges 1$. Hence it is enough to prove that the restriction of \eqref{eq:st} to 
$Y_{ost}$  is the identity. After restricting \eqref{eq:os} to $Y_{ost}$, we deduce 
$$
\tilde{\alpha}_s|_{Y_{ost}}
=\tilde{\alpha}_o|_{Y_{ost}}
=\tilde{\alpha}_t|_{Y_{ost}}
\quad\Rightarrow\quad
{\tilde{\alpha}_t^{-1}\circ\tilde{\alpha}_s|}_{Y_{ost}}=\bone.
$$
Finally, we conclude that the trivialization $\tilde\alpha$ of 
$\pi_*q^*(\eL_\mu^{-1}\otimes\crl V)$ descends to $\cU$, 
as announced. Indeed, we define 
$\bar\alpha:\eO_\cU^{\oplus m}\to\crl V\otimes\eO_\cU$, 
$\bar\alpha(x):=\tilde\alpha_s(x)$ for some $s\in B$ such that  $x\in Y_s$. 
The diagram \eqref{eq:st} implies that $\bar\alpha(x)$ is independent of 
$s\in B$ with $s(x)=0$. 
\end{m-proof}

\begin{m-remark}\label{rmk:glue}
Let $B$ be as in \ref{lm:galois}. The proofs of the Lemmas \ref{lm:pic} and \ref{lm:ind-start} require only the following consequences of \eqref{eq:3nu}: 
\begin{enumerate}
\item[\rm(i)] For all $s\in B$, $\Pic(X)\to\Pic(Y_s)$ is an isomorphism; 
\item[\rm(ii)] For all $s\in B$, there is $B_s\subset B$ dense, such that the intersection $Y_{st}$, $t\in B_s$, is transverse and the diagram \eqref{eq:pic} consists of isomorphisms;
\item[\rm(iii)] $Y_o,Y_{st},Y_{ost}$ are connected, for all $o,s,t\in B$.
\end{enumerate}
This observation will be used in the proof of the Theorem \ref{thm:flag}.
\end{m-remark}

\begin{m-lemma}\label{lm:ind}
Let the situation be as in the Lemma \ref{lm:ind-start}. Then $\crl V$ 
is obtained as a successive extension of line bundles on $X$. 
\end{m-lemma}

\begin{m-proof}
Let ${(q^*\crl V)}_{B}
=\hbox{$\underset{j\in J}\bigoplus$}\,q^*\eL_j\otimes\mbb C^{m_j}$ be an isotypical 
decomposition, with $\eL_j\in\Pic(X)$. 
First we prove the lemma over $\cU$, by induction on the cardinality of $J$. For $|J|=1$, 
we have ${(q^*(\eL^{-1}\otimes\crl V))}_B\cong\eO_{\cY_{B}}^{\oplus m}$ for some 
$\eL\in\Pic(X)$. The Lemma \ref{lm:ind-start} implies   
$\crl V_\cU\cong\eL\otimes\eO_\cU^{\oplus m}$. 

Now suppose that the lemma holds for $|J|\les n$, and prove it for $|J|=n+1$. 
For the maximal elements $M\subset J$, there is a pointwise injective homomorphism 
$
\underset{\mu\in M}{\mbox{$\bigoplus$}}
\eL_\mu\otimes\eO_\cU^{\oplus m_\mu}
\to
\crl V_\cU. 
$ 
Its cokernel $\crl W_\cU$ is locally free over $\cU$ and 
$q^*\crl W_\cU\cong\underset{j\in J\sm M}{\bigoplus}q^*\eL_j^{\oplus m_j}$. 
By the induction hypothesis, $\crl W_\cU$ is obtained by successive extensions 
of $\eL_j$, $j\in J\,\sm\, M$, so the same holds for $\crl V_\cU$. 

It remains to prove that $\crl V$ itself is a successive extension of line bundles on $X$. 
This follows by repeatedly applying the Corollary \ref{cor:ext}. Indeed, each of the successive 
extensions involved in $\crl V_\cU$ uniquely extends to the whole $X$, as $\cU$ is an 
open neighbourhood of $Y_o$, $o\in B$. 
But $\crl V_\cU$ is the result of this process over $\cU$, and $\crl V$ is already defined 
on the whole $X$, so the uniqueness part of the Corollary \ref{cor:ext} yields the conclusion. 
\end{m-proof}

\begin{m-theorem}\label{thm:generic-s}
Let $X$ be an irreducible, smooth, projective $\bk$-variety, and $\eN$ a globally generated, 
ample vector bundle on it satisfying \eqref{eq:3nu}. 
We assume that the restriction of $\crl V$ to the geometric generic section $\mbb Y$ of $\eN$ splits. 
Then the following statements hold:
\begin{enumerate}
\item[\rm(i)] 
If $\bk$ is uncountable, $\crl V$ is a successive extension of line bundles on $X$. 
\item[\rm(ii)] 
Assume that $\bk$ is arbitrary and either one of the following two conditions is fulfilled: 
\begin{enumerate}
\item[\rm(H1)] $H^1(X,\eL)=0$,  for all $\eL\in\Pic(X)$;
\item[\rm(SS)\kern1pt] $\Gamma\bigl(\mbb Y,\cEnd(\crl V_{\mbb Y})\bigr)$ 
is a semi-simple, finite dimensional algebra. 
\end{enumerate}
\nit Then $\crl V$ is a split vector bundle on $X$. 
\end{enumerate}
\end{m-theorem}

\begin{m-proof} 
The proof is done in two steps.\smallskip 

\nit{\it Case $\bk=\mbb C$.} \quad 
Let $B\subset\mbb P\bigl(\Gamma(X,\eN)\bigr)$ be as in the Lemma \ref{lm:ind-start}, 
and decompose 
${(q^*\crl V)}_{B}=\underset{j\in J}{\bigoplus}\,q^*\eL_j\otimes\mbb C^{m_j}$. 
The Lemma \ref{lm:ind} says that $\crl V$ is a successive extension of $\eL_j$, $j\in J$. 
Now assume that either (H1) or (SS) is satisfied. On one hand, 
if $H^1(X,\eL)=0$ for all $\eL\in\Pic(X)$, then any extension of line bundles is trivial, 
so $\crl V$ is isomorphic to 
$\underset{j\in J}{\mbox{$\bigoplus$}}\,\eL_j^{\oplus m_j}$. On the other hand, 
$\Gamma(\mbb Y,\crl E_{\mbb Y})$ is semi-simple if and only if  
$
\Gamma(\mbb Y,\ell_i^{-1}\ell_j)=0
\;\Leftrightarrow\;
\Gamma(Y_s,\eL_i^{-1}\eL_j)=0
$, 
$\forall\,i\neq j\;\forall\,s\in B.$ 
In this case all the elements of $J$ are maximal with respect to \eqref{eq:order}, 
and the conclusion follows from the Lemma \ref{lm:ind-start}. 
\smallskip 

\nit{\it Case $\bk$ arbitrary.}\quad 
Let $\bk_0\subset\bk\cap\mbb C$ be a countable, algebraically closed field, such 
that $X,\eN,\crl V$ are defined over $\bk_0$, and let $X_0,\eN_0,\crl V_0$ be the 
corresponding objects. Then the geometric generic fibres fit into the Cartesian diagram 
$$
\xymatrix@R=1.5em{
\cY_{\bar{\sbbk}}\ar[d]\ar[r]^-\psi&\cY_{\bar{\sbbk}_0}\ar[d]^-\pi
\\
\Spec(\bar{\bbk})\ar[r]&\Spec(\bar{\bbk}_0),
}
$$
and $(q^*\crl V)_{\bar{\sbbk}_0}$ splits by the Proposition \ref{prop:go-down}, so 
$(q^*\crl V)_{\bar{\sbbk}_0}\times_{\bk_0}\mbb C$ splits too. But this latter is 
the restriction of $\crl V_{\mbb C}:=\crl V\times_{\bk_0}{\mbb C}$ to the zero locus of the 
geometric generic section of $\eN_{\mbb C}$, 
hence $\crl V_{\mbb C}$ is a successive extension of line bundles 
$\eL_j$ on $X_\mbb C$, $j\in J$, by the previous step. 

There is an intermediate field $\bk_0\hra\bk_1\hra\mbb C$ of finite type over 
$\bk_0$, such that $\eL_j$, $j\in J$, are defined over $\bk_1$, thus 
$\crl V_0\times_{\bk_0}\bk_1$ is a successive extension of line bundles on 
$X_0\times_{\bk_0}\bk_1$. 

On one hand, if $\bk$ is uncountable, the transcendence degree of $\bk$ over $\bk_0$ 
is infinite because $\bk_0$ is countable. Hence we can realize $\bk_1$ as a sub-field of 
$\bk$, and the conclusion follows by base change. 

On the other hand, if either (H1) or (SS) is fulfilled (over $\bk$), then 
the same holds over $\bk_0$ and $\mbb C$, so $\crl V_{\mbb C}$ splits.  
By applying \ref{prop:go-down} once more, we deduce the splitting of 
$\crl V_0$ on $X_0$ and of $\crl V$ on $X$. 
\end{m-proof}

\begin{m-remark} 
{\rm(i)} 
If $\Gamma(\crl E_{\mbb Y})$ is not semi-simple, then $\crl V$ is a successive extension 
of line bundles on $X$, and we don't know whether $\crl V$ actually splits. The difficulty 
is that the unipotent automorphisms of $\crl V_{\mbb Y}$ act non-trivially on the 
isotypical decompositions of $\crl V_{\mbb Y}$. 

{\rm(ii)} 
We cannot decide the optimality of the factor $1/3$ in \eqref{eq:3nu}. 
The following example illustrates why the triple intersections $Y_{ost}$, 
$o,s,t\in\Gamma(X,\eN)$, are assumed non-empty and connected. 
Let $\crl V=\crl T_{\mbb{P}^2_{\mbb C}}$ be the tangent bundle of 
$X=\mbb{P}^2_{\mbb C}$. It is a non-split, uniform vector bundle of rank two, and its 
restriction to any line $Y\subset\mbb{P}^2_{\mbb C}$ is isomorphic to 
$\eO_Y(2)\oplus\eO_Y(1)$. The incidence variety $\cY$ is the variety of full flags 
in $\mbb C^3$,  and we have the diagram 
$$
\xymatrix@R=1.25em{
\cY\ar[r]^-q\ar[d]_-\pi^-{\mbb{P}^1_{\mbb C}\rm -fibration}&\mbb{P}^2_{\mbb C}
\\
|\eO_{\mbb{P}^2_{\mbb C}}(1)|\cong\mbb{P}^2_{\mbb C} 
}
$$
The geometric generic fibre $\mbb Y$ of $\pi$ is isomorphic to the projective line 
defined over the algebraic closure of the quotient field of $\mbb{P}^2_{\mbb C}$, 
so $q^*\crl T_{\mbb{P}^2_{\mbb C}}$ splits on $\mbb Y$, and there is a ball 
$B\subset|\eO_{\mbb{P}^2_{\mbb C}}(1)|$ such that 
$(q^*\crl T_{\mbb{P}^2_{\mbb C}})|_{\pi^{-1}(B)}$ splits. However, this splitting does 
not descend to $q\bigl(\pi^{-1}(B)\bigr)\subset\mbb{P}^2_{\mbb C}$, for no such $B$. 
Otherwise, the Proposition \ref{prop:formal} would imply that 
$\crl T_{\mbb{P}^2_{\mbb C}}$ splits, a contradiction. 
\end{m-remark}

For {\em arbitrary} varieties defined over 
{\em uncountable} ground fields ({\it e.g.} $\mbb C$), it is enough to check 
the splitting of the restriction to a {\em single} sufficiently general ample subvariety. 

\begin{m-theorem}\label{thm:generic-s2}
Assume that \eqref{eq:3nu} is satisfied, and $\bk$ is {\em uncountable}. 
Let $\bk_0\subset\bk$ be a {\em countable}, algebraically closed sub-field such that 
$X, \eN, \crl V$ are defined over $\bk_0$. 
Consider a regular section $s\in\Gamma(X,\eN)$ with the following properties: 
\begin{itemize}
\item[--] $\crl V_{Y_s}$ is split;
\item[--] in some affine chart induced from $\mbb P_{\bk_0}^N$, the coordinates 
of $[s]\in\mbb P_{\bk}^N$ are algebraically independent over $\bk_0$. 
\end{itemize}
Assume furthermore that either one of the following two conditions is satisfied: 
\\ \null\quad{\rm(H1)}
$H^1(X,\eL)=0\text{ for all }\eL\in\Pic(X)$; 
\\ \null\quad{\rm(SS)}
$\Gamma\bigl(Y_s,\cEnd(\crl V_{Y_s})\bigr)$ is a semi-simple, finite dimensional algebra.\\ 
\nit Then the vector bundle $\crl V$ splits into a direct sum of line bundles on $X$. 
\end{m-theorem}

As $\bk_0$ is countable, the points $[s]\in\mbb P^N_\bk$ with the previous properties 
lie in the complement of a countable union of proper subvarieties of $\mbb P^N_\bk$, 
so we can  reformulate as follows:\smallskip 

{\it If $\bk$ is uncountable, and the restriction to the zero locus of a {\em very general} 
section of $\eN$ splits, then the vector bundle $\crl V$ does the same.}

\begin{m-proof}
Let $(c_1,\ldots,c_N)$ be the coordinates of $[s]$ in the affine chart $c_0\neq 0$ 
on $\mbb P^N_\bk$. By assumption, they are algebraically independent over $\bk_0$, which implies 
$$
\bk_0\subset \bbk_0:=\bk_0(\xi_1,\ldots,\xi_N)
\cong\bk_0(c_1,\ldots,c_N)\subset\bk
\quad\srel{\bk\text{ \small alg. closed}}{\Longrightarrow}\quad 
\bar{\bbk}_0\subset\bk.
$$
(Here $\xi_1,\ldots,\xi_N$ are indeterminates.) 
Therefore the closed point $[s]\in\mbb P^N_\bk$ maps to the generic point of $\mbb P^N_{\bk_0}$. 
Moreover, as $\bk_0$ is countable, it can be realized as a sub-field of $\mbb C$. 
We consider the following diagram: 
$$\xymatrix@R=.55em{
&X\ar[rr]^-b\ar'[d][dd]
&&X_{\bar{\sbbk}_0}\ar[rr]\ar'[d][dd]
&&X_{\bk_0}\ar[dd]
\\ 
Y_s\ar[rr]\ar[dd]_<(.6)\pi\ar[ur]^<(.35)q\ar@{.>}[dr]
&&\cY_{\bar{\sbbk}_0}\ar[dd]_<(.6){\ovl\pi_0}\ar[ur]^<(.2){\ovl q_0}
\ar[rr]\ar@{.>}[dr]
&&\cY_{\bk_0}\ar[dd]_<(.6){\pi_0}\ar[ur]^<(.2){q_0}&
\\ 
&\Spec(\bk)\ar'[r][rr]
&&\Spec(\bar{\bbk}_0)\ar'[r][rr]
&&\Spec(\bk_0)
\\  
[s]\ar[rr]\ar[ur]
&&\mbb P^N_{\bar{\sbbk}_0}\ar[rr]\ar[ur]
&&\mbb P^N_{\bk_0}\ar[ur]
}$$
Now we focus on the diagonal Cartesian rectangle with dotted sides. Our hypothesis 
is that $\crl V_{Y_s}$ splits. Since both 
$\bar{\bbk}_0$ and $\bk$ are algebraically closed, the Proposition \ref{prop:go-down} implies 
that ${(q_0^*\crl V)}_{\bar{\sbbk}_0}$ on $\cY_{\bar{\sbbk}_0}$ splits, 
so ${(q_0^*\crl V)}_{\bar{\sbbk}_0}\times_{\bk_0}\mbb C$ splits too. 
The Theorem \ref{thm:generic-s} implies that $\crl V_\mbb C$ splits 
and we conclude that the initial $\crl V$ is split, by the Proposition \ref{prop:go-down}. 
\end{m-proof}

The theorem can be used to create a variety of applications. It is not clear how to handle 
the following examples by using different methods. 

\begin{m-example}\label{expl:2n} 
(i) (Compare with \ref{expl:grass}(i)) 
Consider $X\srel{\iota}{\hra}\Grs(n;\mbb C^{n+\nu})$ of codimension $c$, 
with $n\ges 4,\nu\ges2, c\les(n-3)\nu-1$, and let $\eG$ be the universal quotient. 
If $X$ satisfies the condition (H1), then the splitting of $\crl V$ along the zero locus 
of a very general section of $\eN:=\eG(1)_X$ implies its global splitting. 
\\[.5ex]
\nit(ii) Let $\crl V$ be a vector bundle on 
$X:=\mbb{P}^m_{\mbb C}\times\mbb{P}^n_{\mbb C}$, 
with $n\ges 2m+1$ and $m\ges 2$. We consider the ample vector bundle 
$\eN:=\crl T_{\mbb{P}^m_{\mbb C}}\boxtimes\eO_{\mbb{P}^n_{\mbb C}}(1)$, 
and a very general section $s\in\Gamma(X,\eN)$. Then $\crl V$ splits if and only if 
its restriction to $Y_s$ does so. (Note that $Y_s\subset X$ has codimension $m$, and 
$Y_s\to\mbb{P}^n_{\mbb C}$ is a $(m+1)$-sheeted ramified covering.) 
\end{m-example}


\section{Splitting along divisors}\label{sct:div}

Here we elaborate on the case where $\eN$ is an ample line bundle and $Y$ is a divisor. 
Throughout $\dim_\bk X\ges 3$, and $\eO_X(1)$ is an ample line bundle on $X$. 

\begin{m-lemma}\label{lm:h1}
Let $D\in|\eO_X(m)|$, with $m\ges 1$, be a divisor such that 
\begin{equation}\label{eq:h1}
H^1(D,\crl E_D(-a))=0,\quad\forall a\ges c.
\quad
(\text {Recall that }\crl E=\cEnd(\crl V).)
\end{equation}
Then hold: 

\nit{\rm(i)} 
The cohomology group $H^1\bigl(X,\crl E(-a)\bigr)$ vanishes for all $a\ges c$. 

\nit{\rm(ii)} 
Assume moreover that $m\ges c$ and $\crl V_D$ splits. Then $\crl V$ splits too. 
\end{m-lemma}

\begin{m-proof}
(i) The Serre vanishing implies: 
$$
a_0{:=}\min\{a\ges c\,|\,H^1(X,\crl E(-j)){=}\,0,\,\forall j\ges a\}<\infty.
$$
If $a_0\ges c+1$, the exact sequence 
$0\,{\to}\eO_X(-m){\to}\,\eO_X{\to}\,\eO_D{\to}\,0$ yields 
$$
\ldots\!\to 
H^1\bigl(\crl E(-m-a_0+1)\bigr)
\to 
H^1\bigl(\crl E(-a_0+1)\bigr)\to
H^1\bigl(\crl E_D(-a_0+1)\bigr)
\to\!\ldots,
$$
with $-m-a_0+1\les-a_0$, $a_0-1\ges c$; the first and last terms vanish, 
so the middle term vanishes too, which contradicts the minimality of $a_0$. 

\nit (ii) As $m\ges c$, the first step implies that 
$\res_D:\Gamma(X,\crl E)\to\Gamma(D,\crl E_D)$ is surjective. 
Hence $\crl V$ splits, by the Lemma \ref{lm:eigv}(ii). 
\end{m-proof}

\begin{m-theorem}\label{thm:gen1}
Assume that either 

\nit{\phantom{ooo}\rm(i)} $D\in |\eO_X(m)|$ is normal, and $\crl E_D(m)$ is ample 
(\emph{e.g.} $\crl E(m)$ is ample), 

\nit{or \rm(ii)} $H^1(D,\ell)=0,\;\forall\,\ell\in\Pic(D)$. 
(See \cite[Proposition 4.13, Corollary 4.14]{ba}.)

\nit Then $\crl V$ splits if and only if $\crl V_D$ splits. 
\end{m-theorem}
The criterion implies that $\crl V$ splits if and only if its 
restriction to a complete intersection surface in $X$ of sufficiently high degree splits. 

\begin{m-proof}
(i) By hypothesis $\crl V_D=\ouset{j=1}{r}{\oplus}\ell_j$ with $\ell_j\in\Pic(D)$. As $\crl E_D(m)$ is ample, $\ell_i^{-1}\ell_j\otimes\eO_D(m+a)$ is ample, for all $i,j$ and $a\ges 0$. The Kodaira vanishing theorem \cite{mum} yields $H^1(\crl E_D(-m-a))=0, \forall a\ges 0$, which is the condition \eqref{eq:h1}. 

\nit(ii) Since $\crl V_D$ splits, $H^1(\crl E_D(-a))=0$ for all $a\ges 0$; 
we conclude as before.
\end{m-proof}

Varieties enjoying additional properties admit stronger splitting criteria. 
\begin{m-definition}\label{horrocks} 
Let $X$ be a scheme and $h\ges 1$ be an integer. We say that $X$ is an 
{\em $h$-splitting scheme} if $H^1(X,\eL)=\ldots=H^{h}(X,\eL)=0$  
for all line bundles $\eL\to X$. 
The cases $h=1,2$ respectively correspond to the notions of {\it splitting} 
and {\it Horrocks scheme} in \cite{ba}. 
\end{m-definition}

If $(X,\eO_X(1))$ is $d$-dimensional, arithmetically Cohen-Macaulay, with 
$\Pic(X)=\mbb Z\cdot\eO_X(1)$, then $X$ is $(d-1)$-splitting. (This follows directly 
from the definition and the Kodaira vanishing.) Examples include Fano varieties with 
cyclic Picard groups (\emph{e.g.} homogeneous spaces $G/P$, with $P\subset G$ a 
maximal parabolic subgroup), and smooth (resp.  very general) complete intersections 
of dimension $d\ges 4$ (resp. $d\ges 3$) in them. 

The next result generalizes \cite[Corollary 4.14]{ba} because we allow $1$- rather than $2$-splitting varieties. 

\begin{m-theorem}\label{thm:gen2} 
Assume $\bk$ is uncountable, $\eO_X(m)$ is globally generated; take $D\in|\eO_X(m)|$ 
very general. (So $D$ is smooth). In either one of the following situations, $\crl V$ splits 
if and only if its restriction $\crl V_D$ splits: 

\nit{\rm(i)} 
$X$ is $2$-splitting, $\dim_\bk(X)=3$, and $\kappa_X(m)$ 
is generated by global sections. (Here $\kappa_X$ stands for the canonical line bundle.)

\nit{\rm(ii)} 
$X$ is $1$-splitting, $\dim_\bk(X)=4$, and $\kappa_X(2m)$ is 
generated by global sections. 

\nit{\rm(iii)} 
$X$ is $1$-splitting, $\dim_\bk(X)\ges5$.
\end{m-theorem}

\begin{m-proof} 
(i) 
The Noether-Lefschetz theorem \cite{rs} states that $\Pic(X)\to\Pic(D)$ is an isomorphism. 
Thus for any $\ell\in\Pic(D)$ there is $\eL\in\Pic(X)$ such that $\eL_D=\ell$. 
The long exact sequence in cohomology associated to  
$0\,{\to}\,\eL(-m){\to}\,\eL{\to}\,\ell{\to}\,0$ yields $H^1(D,\ell)=0$. 
Hence $H^1(D,\crl E_D(-a))=0,\;\forall\,a\ges 0$ because $\crl E_D$ splits, by hypothesis. 
Now apply the Lemma \ref{lm:h1}(ii).

\nit(ii), (iii) The statements are particular cases of the Theorem \ref{thm:generic-s2}. 
\end{m-proof}

An interesting application for vector bundles on the Grassmannian, which combines our results 
obtained so far, is given in \ref{rmk:q43}.

\begin{m-remark}\label{rmk:+char}
The same arguments work in positive characteristics. Assume that $\chr(\bk)>\dim_\bk X\ges 4$, 
the pair $(X,\eO_X(1))$ admits a $W_2(\bk)$-lifting, and $X$ is $2$-splitting. Then $\crl V$ splits 
if and only if $\crl V_D$ splits. 
\end{m-remark}


\section{Splitting of vector bundles on Grassmannians \break 
(The ampleness of $\eN$ is necessary?)}\label{sct:grass}

Throughout this article we restricted $\crl V$ to `test subvarieties' which are zero loci of regular sections in an \emph{ample} vector bundle $\eN$. However, our strategy to prove the splitting of $\crl V$ (cf. page \pageref{strategy}) makes sense without this restriction, and is natural to ask whether the hypothesis can be weakened. In the general setting of the Section \ref{sct:x-xy}, this does not seem possible, because the ampleness hypothesis is used to deduce the vanishing of certain cohomology groups. However, the answer to the question appears to be affirmative if $\eN$ is globally generated and one has \textit{a priori} information about the zero loci of its regular sections. We illustrate our statement with two concrete examples.


\subsection{The case of the Grassmannian}\label{ssct:gr}

\begin{m-theorem}\label{thm:grass}
A vector bundle $\crl V$ on $X:=\Grs(e;\bk^d)$, with $e,d-e\ges2$, is split if and only if its restriction to $Y:=\{U\in\Grs(e;\bk^d)\mid \bk^{e-2}{\subset}\,U{\subset}\,\bk^{e+2}\}\cong\Grs(2;\bk^4)$ is so.
\end{m-theorem}

Cohomological splitting criteria for vector bundles on Grassmannians have been obtained in \cite{ott,mal,ar-ma}. However, they involve many conditions. We believe that our result is interesting for its simplicity: it reduces the problem of splitting on the Grassmannian, which is a high dimensional object, to a $3$-dimensional quadric $Q_3\in\mbb P^4$. (See the Remark \ref{rmk:q43} for the reduction from $\Grs(2;4)$ to $Q_3$.)

For $e=2, d=4$ there is nothing to prove. For $e=2,d\ges5$, we use the duality $\Grs(e;\bk^d)\cong\Grs(d-e;\bk^d)$, so we can write $d=\nu+n+1$ and $e=n+1$, with $\nu,n\ges2$. The theorem is obtained by repeatedly applying the following:

\begin{m-proposition}\label{prop:grass}
A vector bundle $\crl V$ on $X=\Grs(n+1;W)$, with $W\cong\bk^{\nu+n+1}$, 
$\nu,n\ges2$, splits if and only if its restriction to some smaller Grassmannian 
$Y\cong\Grs(n;\nu+n)$ contained in $X$ splits. 
\end{m-proposition}

The universal quotient $W\otimes\eO_X\,{\srel{\beta}{\to}}\,\eN$ induces 
an isomorphism $W\,{\srel{\cong}{\to}}\,\Gamma(\eN)$, and $s\in W\sm\{0\}$ 
determines a section in $\eN$ whose zero locus is the `smaller' Grassmannian 
$Y:=\Grs(n;W_s)$, with $W_s:=W/\lran{s}$. 

We remark that $\eN$ is \emph{not ample} on $X$ because its restriction to any (straight) 
line $l\subset X$ is isomorphic to $\eO_l^{\nu-1}\oplus\eO_l(1)$. 

\begin{m-proof} 
We follow the strategy described on the page \pageref{strategy}. 

\nit{\it Claim~1}\quad  
For any vector bundle $\euf F$ on $X$ holds: 
\begin{equation}\label{eq:cdgrs}\kern-2.25ex
\begin{array}{l}
\bullet\,H^1(X,\euf F\otimes\eI_Y^m)=0,\text{ so }
\Gamma(X,\euf F)\to\Gamma(Y_m,\euf F_{Y_m})\text{ is surjective},\,\forall m\gg0;
\\[1ex] 
\bullet\,\text{the cohomological dimension }\cd(X\sm Y)\les\dim X-(n+1). 
\end{array}\kern-1ex
\end{equation}
Let $\tld X:=\Bl_Y(X)$ be the blow-up of $X$ along $Y$ and $\pi\,{:}\,\tld X{\to}\,X$ 
be the projection; we denote the exceptional divisor by $E\,{=}\,\mbb P(\eN_Y)\,{\subset}\,\tld X$. 
As $Y$ is the zero locus of $s\in\Gamma(\eN)$, we have $\tld X\,{\subset}\,\mbb P(\eN)$, and holds:

\nit{\rm(i)} 
$\eO_{\tld X}(-E)=\eO_{\mbb P(\eN)}(1)|_{\tld X}$ is $\pi$-relatively ample;

\nit{\rm(ii)} 
$H^1(X,\eF\otimes\eI_Y^m)=\,H^1(\tld X,\pi^*\eF\otimes\eO_{\tld X}(-mE))$, 
$\forall m\ges 1$.

\nit Furthermore, as $\eN\cong\overset{\nu-1}{\bigwedge}\eN^\vee\otimes\det(\eN)$ 
and $\eN$ is globally generated, we have 
$$
\mbb P(\eN)\cong\mbb P(\overset{\nu-1}{\bigwedge}\eN^\vee\otimes\det(\eN))
\subset
\mbb P\bigl(\overset{\nu-1}{\bigwedge}W^\vee\otimes\eO_X(1)\bigr)
\cong X\times\mbb P,
$$
where $\mbb P{:=}\,\mbb P(\overset{\nu-1}{\bigwedge}W^\vee)$, and 
$\eO_{\mbb P(\eN)}(1){=}\,
\bigl(\eO_X(-1){\boxtimes}\,\eO_{\mbb P}(1)\bigr)|_{\mbb P(\eN)}$.  

Pointwise, the morphism $\mbb P(\eN)\to\mbb P$ is defined by 
\begin{equation}\label{eq:q}
(x,\lran{e_x})\mt 
\det(\eN_x/\lran{e_x})^\vee\subset
\overset{\nu-1}{\bigwedge}\eN_x^\vee\subset
\overset{\nu-1}{\bigwedge}W^\vee,
\end{equation}
where $\lran{e_x}$ stands for the line generated by $e_x\in\eN_x$. 
Its restriction  $q:\tld X\to\mbb P$ to $\tld X$ corresponds to the commutative diagram 
$$
\xymatrix@R=1.5em{
0\ar[r]&\eO_X\ar[r]^-{s}\ar@{=}[d]&
W\otimes\eO_X\ar[r]\ar@{->>}[d]^-{\;\beta}&W_s\otimes\eO_X\ar[r]\ar@{->>}[d]&0
\\ 
&\eO_X\ar[r]^-{\beta s}&\eN\ar[r]&\eN/\lran{\beta s}\ar[r]&0.
}
$$
The homomorphism $\beta s$ is injective precisely over $X\sm Y$ 
and \eqref{eq:q} shows that $q$ is the desingularization of the rational map 
\begin{equation}\label{eq:GG}
\Grs(n+1;W)\dashto\Grs(n+1;W_s),\quad 
[U\subset W]\mt [(U+\lran{s})/\lran{s}\subset W_s],
\end{equation}
followed by the usual Pl\"ucker embedding of $\Grs(n+1;W_s)$. 

For an arbitrary vector bundle $\tld\eF$ on $\tld X$, we have 
$R^jq_*(\tld\eF\otimes\pi^*\eO_X(m))=0$, for all $j>0$ and $m\gg0$, 
so the Leray spectral sequence yields 
\begin{equation}\label{eq:h}
\begin{array}{rl}
&
H^{\dim X-1}\bigl(\tld X,
\tld\eF\otimes\pi^*\eO_X(m)\otimes q^*\eO_{\mbb P}(-m)
\bigr)
\\ 
=
&
H^{\dim X-1}\bigl(\Grs(n+1;W_s),
q_*(\tld\eF\otimes\pi^*\eO_X(m))\otimes\eO_{\mbb P}(-m)
\bigr)
=0.
\end{array}
\end{equation}
The last equality holds because $\dim X-1>\dim\bigl(\Grs(n+1;W_s)\bigr)$. 
The Serre duality implies 
$H^1(X,\euf F\otimes\eI_Y^m)
\,{=}\,
H^1(\tld X,\pi^*\eF(-m)\otimes q^*\eO_{\mbb P}(m))
\!\srel{\eqref{eq:h}}{=}\!0$, for $m\gg0$. 

It remains to estimate the cohomological dimension of $X\sm Y$. For this, we observe that the morphism $q:X\sm Y\to\Grs(n+1;W_s)$ is affine (cf. \eqref{eq:GG}), so 
$\cd(X\sm Y)\les\dim\Grs(n+1;W_s)=\dim X-(n+1)$.  
\smallskip 

\nit{\it Claim~2}\quad 
For any line bundle $\ell$ on $Y$ and $m\ges 1$ holds 
$$
H^1(Y,\Sym^m(\eN_Y^\vee)\otimes\ell)=0\quad \text{(cf. \ref{lm:split}(ii))}.
$$
Note that $\det(\eN_Y)=\eO_Y(1)$ generates $\Pic(Y)$, so  
$\ell=\eO_Y(k)$ for some $k\in\mbb Z$. The cohomology group above can be computed 
on $F:=\mbb P(\eN_Y^\vee)\srel{f}{\to}Y$: 
\begin{equation}\label{h1}
H^1(Y,\Sym^m\eN_Y^\vee\otimes\eO_Y(k))
\cong 
H^{\nu}(F,\underbrace{\eO_f(-\nu-m)\otimes\eO_Y(k+1)}_{=:\eL}).
\end{equation}
Also, $F$ is isomorphic to the (homogeneous) variety of partial flags 
$$
0\,{\subset}\,U_{n}\,{\subset}\,U_{\nu+n-1}\,{\subset}\,W_s^{\nu+n},
$$ 
where the indices indicate the dimensions. The vanishing 
of \eqref{h1} follows from Bott's theorem \cite{bott,de}. 
We follow \cite[Chapter 4]{we} which treats in detail the flag varieties for the general linear group. 
Let $\eU_{n}, \eU_{\nu+n-1}$ be the tautological bundles on $F$ of indicated ranks. Then 
$$
\eO_f(1)\cong (W_s\otimes\eO_Y)/\eU_{\nu+n-1},\; f^*\eO_Y(1)\cong\det(\eU_{n})^{-1},
$$
so $\eL$ corresponds to the weight 
$a=(\underbrace{k+1,\ldots,k+1}_{n\;\text{times}},
\underbrace{0,\ldots,0}_{\nu-1\text{ times}},\nu+m)\in\mbb Z^{\nu+n}$ 
(cf. \cite[pp. 112]{we}). 
We denote $\rho:=(\nu+n-1,\ldots,0)$. Bott's theorem says that $H^\nu(F,\eL)=0$ 
if and only if either one of the following two cases occur: 
\begin{enumerate}
\item 
$a+\rho$ is singular---that is, it contains two identical entries;
\item 
$a+\rho$ is non-singular and the number of strict order inversions (of the decreasing order) 
in $a+\rho$ is different of $\nu$.
\end{enumerate}
We claim that, for all $m\ges 1$, 
$a+\rho=(\,\nu+k+n,\ldots,\nu+k+1,\nu-1,\ldots,1,\nu+m\,),$ 
fulfils one of the two conditions above. 
Note that $\eL^{-1}$ is ample for $k+1<0$, so \eqref{eq:h1} vanishes by 
Kodaira's theorem. For $k+1\ges 0$, the first $\nu+n-1$ 
terms of $a+\rho$ are strictly decreasing, hence only the last term can contribute to a 
strict inversion. Actually, $\nu+m>\nu-1$, so we do have $\nu-1$ inversions. 
There is exactly one more strict inversion precisely for $\nu+m=(\nu+k+1)+1$, that is 
$m=k+2$. But in this case $a+\rho$ is singular since $n\ges 2$, by hypothesis. 
\smallskip 

\nit{\it Claim~3}\quad 
Assume that $\crl V_Y$ splits. Then 
$\Gamma(Y_m,\crl E_{Y_m})\to \Gamma(Y,\crl E_Y)$ is surjective, $\forall m\ges 1$. 
Indeed, we tensor by $\crl E_Y$ the exact sequence:
$$
0\,{\to}\Sym^m\eN_Y^\vee\cong\eI_Y^{m}/\eI_Y^{m+1}
{\to}\,\eO_{Y_m}{\to}\,\eO_{Y_{m-1}}{\to}\,0.
$$
Since $\crl E_Y$ is a direct sum of line bundles, we apply the Claim~2.
\end{m-proof}

\begin{m-remark}\label{rmk:q43}
The results obtained in the Section \ref{sct:div} allow to probe the splitting of the vector 
bundle $\crl V$ on even lower dimensional subvarieties. Indeed, 
the image of $\Grs(2;4)$ by the Pl\"ucker embedding is the smooth $4$-dimensional 
quadric $Q_4\subset\mbb P^5$. 
Let $Q_3{\subset}\,Q_4$ be an arbitrary smooth hyperplane section, and 
$S{\subset}\,\mbb P^4$ be a very general intersection of $Q_3$ with a 
quartic in $\mbb P^4$. (Observe that $S\subset Q_4$ is a surface with 
$\Pic(S)=\mbb Z\cdot\eO_S(1)$, and $\kappa_S=\eO_S(1)$.) 
The Theorems \ref{thm:gen1}(ii) and \ref{thm:gen2}(i)(a) respectively imply that 
$\crl V$ splits on $Q_4$ if and only if either $\crl V_{Q_3}$ or $\crl V_S$ splits. 
(In the latter case we require $\bk$ to be uncountable.)
\end{m-remark}

\subsection{The case of partial flag varieties}\label{ssct:fl}

Here we deduce a splitting criterion for vector bundles on partial flag varieties. To the author's knowledge, there are no previously known results in this case. First we introduce the notation. Let 
$$
d_\bullet:=(0=d_0< d_1<\dots <d_t<d_{t+1})
$$
be a strictly increasing sequence of integers, and $\nu_j:=d_{t+1}-d_j$, for $j=0,\ldots, t$. 
For shorthand, we denote $2_\bullet:=(\,0<2<\dots<2t<2(t+1)\,)$. 

If $d'_\bullet=(0=d'_0<d'_1<\dots<d'_t<d'_{t+1})$ is another sequence, we write: 
\begin{equation}\label{eq:<}
d'_\bullet\les d_\bullet
\;\;\Leftrightarrow\;\;
d'_j-d'_{j-1}\les d_{j}-d_{j-1},\;\forall 1\les j\les t+1.
\end{equation}
We consider the partial flag variety 
$$
F_{d_\bullet}=\Flag(d_1,...,d_t;d_{t+1}):=\big\{
U_\bullet=(0\subset U_{d_1}\subset\ldots\subset U_{d_t}\subset W:=\bk^{d_{t+1}})
\big\},
$$
where the indices indicate the dimensions of the vector spaces. For $d'_\bullet\les d_\bullet$, a choice of subspaces  $\bk^{d_1-d'_1}\subset\dots\subset\bk^{d_{t+1}-d'_{t+1}}\subset\bk^{d_{t+1}}$ yields the embedding: 
$$
\begin{array}{l}
\iota:F_{d'_\bullet}\to F_{d_\bullet},
\\[1ex] 
U_{d'_\bullet}
\mt 
\big(\bk^{d_1-d'_1}\oplus U_{d'_1}
\subset{\dots}\subset
\bk^{d_t-d'_t}\oplus U_{d'_t}
\subset
\bk^{d_{t+1}-d'_{t+1}}\oplus\bk^{d'_{t+1}}=W
\big).
\end{array}
$$
There are $t$ tautological bundles $\eU_{d_j}$ and  $t$ universal quotient bundles $\eN_j$ on $F_{d_\bullet}$, with $\rk(\eU_{d_j})=d_j$, $\rk(\eN_j)=\nu_j$, for $j=1,\ldots,t$. 

As usual, throughout the section, $\crl V$ is a vector bundle on $F_{d_\bullet}$ and $\crl E:=\cEnd(\crl V)$. 

\begin{m-theorem}\label{thm:flag} 
Assume that the sequence $d_\bullet$ satisfies $d_\bullet\ges 2_\bullet$. Then the vector bundle $\crl V$ on $F_{d_\bullet}$ splits if and only if it does along some $\iota\bigl(F_{2_\bullet}\bigr)$. 
\end{m-theorem} 
Everything is defined over a countable, algebraically closed field, which is simultaneously a subfield of $\bk$ and of $\mbb C$ (cf. Section \ref{sct:gener}). Using twice the invariance of the splitting under base change (cf. Proposition \ref{prop:go-down}), we may---and we henceforth do---assume that $\bk=\mbb C$. 
The proof of the statement requires a few intermediate results. 

\begin{m-lemma}\label{lm:cohom0}
With the previous notation, assume that $d_{j+1}-d_{j-1}\ges 3$, for all $j=1,\ldots,t$ (\textit{e.g.} $d_\bullet\ges 2_\bullet$). Then $F_{d_\bullet}$ is $1$-splitting (cf. Definition \ref{horrocks}). 
\end{m-lemma}

\begin{m-proof}
Any line bundle on $F_{d_\bullet}$ can be written in the form 
$$
\ell=\det(\eU_{d_1})^{-a_1}\otimes\dots\otimes\det(\eO_F^{d_{t+1}}/\eU_{d_t})^{-a_t}, 
\;\text{with}\;a_1,\dots,a_t\in\mbb Z.
$$
It corresponds to the weight 
$
\alpha=\big(
\underbrace{a_1,...,a_1}_{d_1\;\text{times}},
\underbrace{a_2,...,a_2}_{d_2-d_1\;\text{times}},...,\kern-1ex
\underbrace{a_t,...,a_t}_{d_{t+1}-d_t\;\text{times}}
\kern-1.5ex\big)\in\mbb Z^{d_{t+1}}.
$

Let $\rho:=(d_{t+1}-1,\dots,1,0)$. By Bott's theorem, $H^1(F,\ell)\neq0$ if and only if 
$$
\alpha+\rho=\big(
\underbrace{a_1+d_{t+1}-1,\dots,a_1+d_{t+1}-d_1}_{\text{length}=d_1}
,\dots,
\underbrace{a_t+d_{t+1}-d_t-1,\dots,a_t}_{\text{length}=d_{t+1}-d_t}
\big)
$$
is \emph{non-singular} (its entries are pairwise distinct), and it contains \emph{exactly one inversion} (for the decreasing order). The $t$ blocks which compose $\alpha+\rho$ are strictly decreasing and the only way to simultaneously achieve the previous two conditions is to have two consecutive blocks of length one each, that is  $d_{j_0+1}-d_{j_0}=d_{j_0}-d_{j_0-1}=1$ for some $j_0$. 
\end{m-proof}

For $1\les j\les t+1$, let $d^j_\bullet:=(d_1,...,d_{j-1},d_j-1,...,d_t-1;d_{t+1}-1).$ We denote 
$$
S:=\{(e,\eta)\in W\oplus W^\vee\mid\eta(e)\neq 0\}.
$$
An element $s=(e,\eta)\in S$ determines a section in $\eN_j\oplus\eU_{d_{j-1}}^\vee$, which is globally generated, whose zero locus is 
\begin{equation}\label{eq:eeta}
Y_s=Y_{(e,\eta)}=
\{U_\bullet\in F_{d_\bullet}\mid U_{d_{j-1}}\subset\Ker(\eta),\; e\in U_{d_j}\}
\cong F_{d^j_\bullet}.
\end{equation}
Therefore we have the same situation as in the section \ref{sct:gener}: 
\begin{equation}\label{eq:Fdj}
\xymatrix@R=.2em@C=2.5em{
&\cY\ar[dl]_-\pi\ar[dr]^-q&
\\ 
S&&F_{d_\bullet},
}
\end{equation}
where $\cY$ is the zero locus of the universal section in $\eN_j\oplus\eU_{d_{j-1}}^\vee$ on $S\times F_{d_\bullet}$. 


\begin{m-lemma}\label{lm:cd}
Assume that $d_\bullet\ges 2_\bullet$ and there is an index $1\les j\les t+1$ such that $d_j-d_{j-1}\ges 3$. Then the cohomological dimension of $F_{d_\bullet}\sm F_{d^j_\bullet}$ satisfies:  $\,\cd(F_{d_\bullet}\sm F_{d^j_\bullet})\les\dim F_{d_\bullet}-3$. 
\end{m-lemma}

\begin{m-proof}
For $Y_{(e,\eta)}$ as in \eqref{eq:eeta}, we have: 
$$
\begin{array}{ll}
O:=F_{d_\bullet}\sm Y_{(e,\eta)}
&=
\big(\,F_{d_\bullet}\sm
\{\underbrace{U_\bullet\mid e\in U_{d_j}}_{=:F'}\}\,\big)
\cup 
\big(\,F_{d_\bullet}\sm
\{\underbrace{U_\bullet\mid U_{d_{j-1}}\subset\Ker(\eta)}_{=:F''}\}\,\big)
\\ 
&=O'\cup O''.
\end{array}
$$
Note that  $F'=\pr_\Grs^{-1}\big(\Grs(d_j-1;d_{t+1}-1)\big)$, where $\pr_\Grs:F_{d_\bullet}\to\Grs(d_j;d_{t+1})$ is the natural projection. The Leray spectral sequence for $\pr_\Grs$ and \eqref{eq:cdgrs} imply that $\cd(O')=\cd(F_{d_\bullet}\sm F')\les\dim F_{d_\bullet}-d_j.$ This proves the lemma for $j=1$. 

By using the duality $\Flag(d_1,...;d_{t+1})\cong\Flag(d_{t+1}-d_t,...;d_{t+1})$ we deduce  $\cd(O'')\les\dim F_{d_\bullet}-(d_{t+1}-d_{j-1}).$ This proves the lemma for $j=t+1$. 

Now assume $2\les j\les t$. The morphism $\pr:F_{d_\bullet}\srel{}{\to}G:=\Flag(d_{j-1},d_j;d_{t+1})$ is smooth, projective, and $Y_{(e,\eta)}$ is the pre-image of the analogous $Y'_{(e,\eta)}\subset G$. The Leray spectral sequence for $F_{d_\bullet}\sm Y_{(e,\eta)}\to G\sm Y'_{(e,\eta)}$ shows that it is enough to prove that $\cd(G\sm Y'_{(e,\eta)})\les\dim G-3$. 
Henceforth, we assume that
$$
d_\bullet=(0<a<b<d),\text{ with }a\ges 2,\,b-a\ges3,\,d-b\ges2.
$$ 
 
For a coherent sheaf $\eG$ on $O$, we have the Mayer-Vietoris sequence 
$$
\ldots\to H^{i-1}(O'\cap O'',\eG)\to H^{i}(O,\eG)\to H^{i}(O',\eG)\oplus H^{i}(O'',\eG)\to\ldots
$$ 
The previous discussion shows that $\cd(O'),\,\cd(O'')\les\dim F_{d_\bullet}-3$, so it is enough to prove $\cd(O'\cap O'')\les\dim F_{d_\bullet}-4$. 

The pair $(e,\eta)$ decomposes $W=\Ker(\eta)\oplus\lran{e}$; let $\pi:W\srel{}{\to}\Ker(\eta)$ be the projection. Since 
$O'\cap O''=\{U_\bullet\mid U_a\not\subset\Ker(\eta),\;e\not\in U_b\},$ 
we deduce the morphism:  
$$
f:O'\cap O''\to G':=\Flag\big(a-1,b;\Ker(\eta)\big),\quad
[U_a\subset U_b]\mt[U_a\cap\Ker(\eta)\subset\pi(U_b)].
$$
\nit\textit{Claim~1}\quad The fibre over $[V_{a-1}\subset V_b]\in G'$ is isomorphic to $\big(\frac{V_b}{V_{a-1}}\big)^\vee\sm\{0\}$; its cohomological dimension equals $b-a$. 

Indeed, let $[U_a\subset U_b]$ be in the fibre. Since $\pi:U_b\to V_b$ is an isomorphism, $U_b$ is the graph of a (uniquely defined) homomorphism $h:V_b\to\lran{e}$: 
$$
U_b=\{(v,h(v))\mid v\in V_b\}.
$$
Also, we have $U_a=V_{a-1}+\mbb C\cdot(v_0,e)$, with $v_0\in\pi(U_a)\subset\pi(U_b)=V_b$. The inclusion $U_a\subset U_b$ implies $h(V_{a-1})=0$, $h(v_0)=e$, hence $h\in\big(\frac{V_b}{V_{a-1}}\big)^\vee\sm\{0\}$. Conversely, any such $h$ defines a flag $[U_a\subset U_b]\in O'\cap O''$. \smallskip

\nit\textit{Claim~2}\quad $\cd(O'\cap O'')\les\dim F_{d_\bullet}-4$. 

For a quasi-coherent sheaf $\eG$ on $O'\cap O''$, the previous step implies that $R^{>(b-a)}f_*\eG=0$. Leray's spectral sequence yields: $\cd(O'\cap O'')\les (b-a)+\dim G'=\dim F_{d_\bullet}-(b-a+1)$. 
\end{m-proof}


\begin{m-lemma}\label{lm:kod-spen}
Assume that $d_\bullet\ges 2_\bullet$ and $d_j-d_{j-1}\ges 3$ for some $1\les j\les t+1$. We assume that there is $o\in S$ such that $\crl V_{Y_o}$ splits. Then there is an open ball $B\subset S$, such that $(q^*\crl V)_B$ splits.
\end{m-lemma}

\begin{m-proof}
The claim follows from the Kodaira-Spencer deformation theory. 
Since $d^j_\bullet\ges 2_\bullet$, the Lemma \ref{lm:cohom0} implies that $Y_o$ is $1$-splitting. 
The fibration $\cY\to S$ is locally trivial (in the analytic topology) with fibres isomorphic to $F_{d^j_\bullet}$, hence $q^*\crl V$ is locally an analytic family of vector bundles on $Y_o\cong F_{d^j_\bullet}$. 

Since $\crl E_{Y_o}$ splits, $H^1(Y_o,\crl E_{Y_o})=0$, therefore $\crl V_{Y_o}$ is rigid. Hence there is a ball $B\subset S$, such that $\crl V_{Y_s}\cong\crl V_{Y_o}$, for $s\in B$ (cf. \cite[Theorem 7.4]{ks}, \cite[Theorem 2.7]{hart-deform}). Possibly after shrinking $B$, the splitting of $\crl V_{Y_o}$ extends over $\cY_B$. 
\end{m-proof}

\begin{m-proof}(of the Theorem \ref{thm:flag})\quad 
We prove by induction on $t+\ouset{j=1}{t+1}{\sum}d_j$. 
For $t=1$ and $\rk(\crl V)$ arbitrary, it is the Theorem \ref{thm:grass}. 
For the inductive step, let $j$ be minimal such that $d_j-d_{j-1}\ges3$. 
Since $F_{2_\bullet}\subset F_{d^j_\bullet}$, the induction hypothesis implies that $\crl V$ splits along $Y_o=\iota(F_{d^j_\bullet})\subset F_{d_\bullet}$. By the Lemma \ref{lm:kod-spen}, there is an open ball $o\in B\subset S$ such that $(q^*\crl V)_B$ splits. 

\nit\textit{Claim}\quad 
The gluing procedure in the Lemma \ref{lm:ind-start} applies, so there is an open analytic neighbourhood $\cU$ of $Y_o$ such that $\crl V_\cU$ is a successive extension of line bundles on $\cU$. 
Indeed, as we pointed out in the Remark \ref{rmk:glue}, the proof of the Lemma \ref{lm:ind-start} is based on the following assumptions: 
\begin{enumerate}
\item[(i)] 
the diagram \eqref{eq:pic}---which compares the Picard groups of $F_{d_\bullet},Y_s,Y_{st}=Y_s\cap Y_t$, for general $s,t\in S$---should consist of isomorphisms;
\item[(ii)] 
the triple intersections $Y_{ost}=Y_o\cap Y_s\cap Y_t$ should be connected. 
\end{enumerate}
We verify that these conditions are fulfilled. For (i), take $s=(e,\eta), s'=(e',\eta')$ generic, so the intersection $Y_{st}$ is transverse, and observe that: 
$$
\begin{array}{ll}
Y_{st}&
=\{\,U_\bullet\mid e,e'\in U_{d_j}\text{ and }U_{d_{j-1}}\subset\Ker(\eta)\cap\Ker(\eta')\,\}
\\[1ex]&
\cong\Flag(d_1,...,d_{j-1},d_j-2,...;d_{t+1}-2).
\end{array}
$$
Clearly, the restrictions $\Pic(F_{d_\bullet})\to\Pic(Y_s)\to\Pic(Y_{st})$ are isomorphisms. For (ii), note that $Y_{ost}\cong\Flag(d_1,...,d_{j-1},d_j-3,...;d_{t+1}-3)$ is connected, even for $d_{j-1}=d_j-3$. 

Let $\hat F$ be the formal completion of $F_{d_\bullet}$ along $Y_o$. The claim implies that $\crl V_{\hat F}$ is a successive extension of line bundles. But $\cd(F_{d_\bullet}\sm Y_o)\les\dim F_{d_\bullet}-3$, by the Lemma \ref{lm:cd}, so the following are isomorphisms, cf. \cite[Theorem III.3.4(b)]{ha}:
$$
\,0=\Ext^1(\ell,\ell')\srel{\cong}{\to}\Ext^1(\ell_{\hat F},\ell'_{\hat F}),
\,\forall\ell,\ell'\in\Pic(F_{d_\bullet}),
\;
\Gamma(F_{d_\bullet},\crl E)\srel{\cong}{\to}\Gamma(\hat F,\crl E_{\hat F}).
$$ 
We deduce that $\crl V_{\hat F}$ splits and consequently $\crl V$ splits too, by the Lemma \ref{lm:eigv}.
\end{m-proof}


\end{document}